\documentclass{siamart0516}
\usepackage[margin=1in]{geometry}
\usepackage{amsmath,booktabs}
\usepackage{amssymb}
\usepackage{amsfonts} %for creating real number symbol. 
\usepackage{todonotes}
\usepackage[normal]{subfigure}
\usepackage{color}

\providecommand{\U}[1]{\protect\rule{.1in}{.1in}}

\renewcommand\thesection{\arabic{section}}
\renewcommand\thesubsection{\thesection.\arabic{subsection}}
\renewcommand\thesubsubsection{\thesubsection.\arabic{subsubsection}}
\renewcommand\theequation{\thesection.\arabic{equation}}
\renewcommand{\eqref}[1]{(\ref{#1})}
\usepackage[numbers,sort&compress]{natbib}

\newcommand{\eps}{{\displaystyle \varepsilon}}
\newcommand{\bsub}{\begin{subequations}}
\newcommand{\esub}{\end{subequations}$\!$}
\newcommand{\ds}[0]{\displaystyle}
\newcommand{\bs}[0]{\boldsymbol}
\newcommand{\littleoh}{ \mbox{{\scriptsize $\mathcal{O}$}}}
\newcommand{\bigoh}{\mathcal{O}}

\newcommand{\bx}{\mathbf{x}}

\newcommand{\by}{\mathbf{y}}

\newcommand{\bg}{\mathbf{g}}

\newcommand{\bxi}{\boldsymbol{\xi}}

\newcommand{\Om}{\Omega}
\newcommand{\pOm} {\partial \Omega}

\newcommand{\A}{{\mathcal{A}}}

\newcommand\RR{{\mathbb R}}
\newcommand\CC{{\mathbb C}}

\newcommand{\bd}{{\partial}}
\newcommand{\DD}{{\mathcal{D}}}
\newcommand{\JJ}{{\mathcal{J}}}
\newcommand{\nn}{{\mathbf{n}}}
\newcommand{\pderiv}[2]{\frac{\partial #1}{\partial #2}}
\newcommand{\xx}{{\mathbf{x}}}
\newcommand{\yy}{{\mathbf{y}}}

\begin{document}

\title{Trapping of planar Brownian motion: Full first passage time distributions by Kinetic Monte-Carlo, asymptotic and boundary integral methods.}

\date{\today} 
\author{Jake Cherry, Alan E. Lindsay, Adri\'{a}n Navarro Hern\'{a}ndez and Bryan Quaife}

\baselineskip=16pt

\maketitle

\begin{abstract}
  We consider the problem of determining the arrival statistics of
  unbiased planar random walkers to complex target configurations. In contrast to problems posed in finite domains, simple
  moments of the distribution, such as the mean (MFPT) and variance, are
  not defined and it is necessary to obtain the full arrival statistics.
  We describe several methods to obtain these distributions and other associated quantities such as splitting probabilities. One approach combines a
  Laplace transform of the underlying parabolic equation with matched
  asymptotic analysis followed by numerical transform inversion. The
  second approach is similar, but uses a boundary integral equation
  method to solve for the Laplace transformed variable. To validate the
  results of this theory, and to obtain the arrival time statistics in
  very general configurations of absorbers, we introduce an efficient
  Kinetic Monte Carlo (KMC) method that describes trajectories as a combination
  of large but exactly solvable projection steps. The effectiveness of
  these methodologies is demonstrated on a variety of challenging
  examples highlighting the applicability of these methods to a variety
  of practical scenarios, such as source inference. A particularly useful finding arising from
  these results is that homogenization theories, in which complex
  configurations are replaced by equivalent simple ones, are remarkably
  effective at describing arrival time statistics.
\end{abstract}

%%%%%%%%%%%%%%%%%%%%%%%%%%%%%%%%%%%%%%%%%%%%%%%%%%%%%%%%%%%%%%%%%%%%%%%%

\label{firstpage}



\pagestyle{myheadings}
\markboth{Jake Cherry, Alan E. Lindsay, Adri\'{a}n Navarro Hern\'{a}ndez and Bryan Quaife}{Arrival times of planar diffusion}

\section{Introduction}\label{sec:intro}

\begin{figure}[htbp]
\centering
\includegraphics[width=0.45\textwidth]{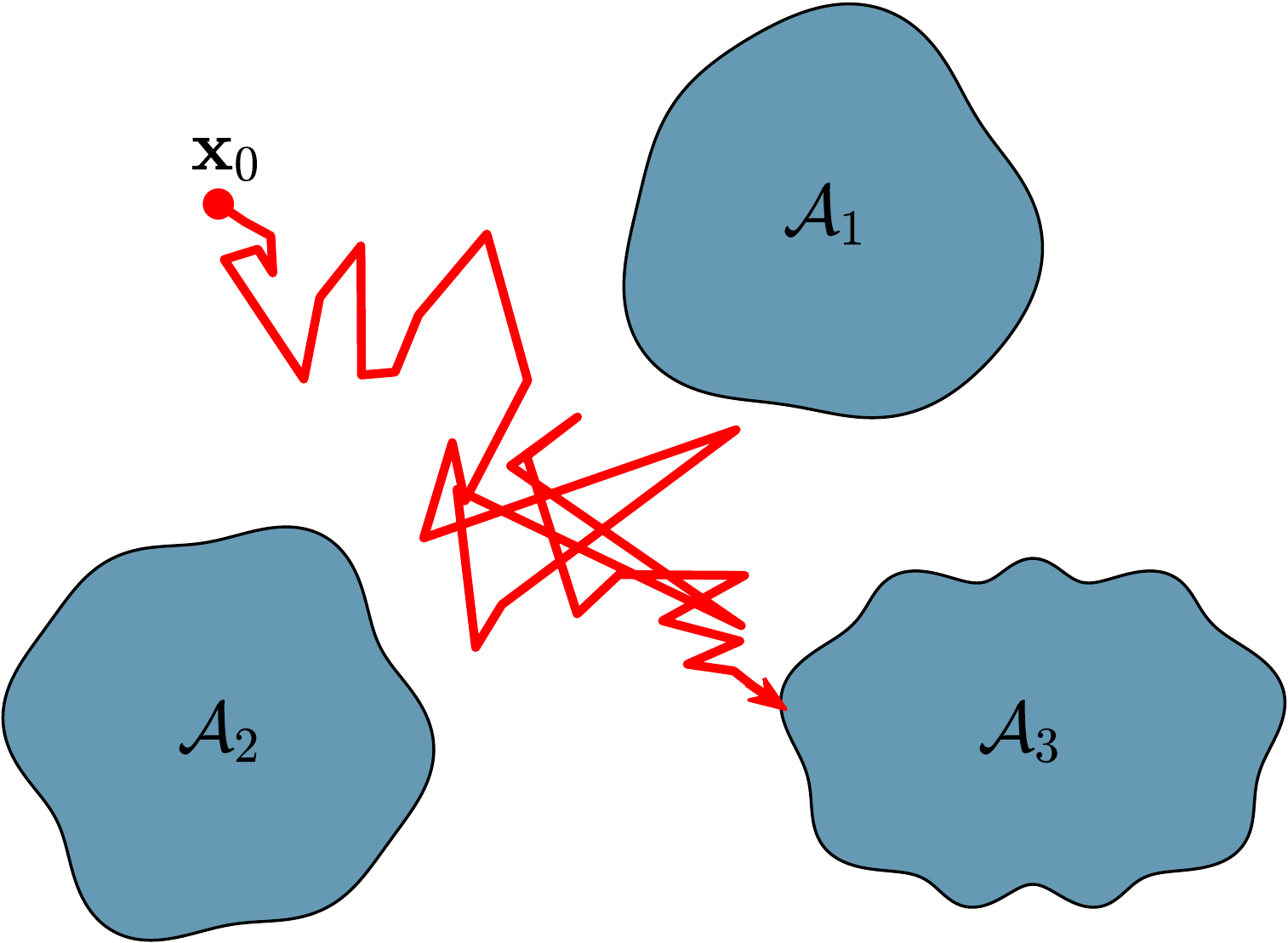}
\caption{Schematic of planar Brownian motion to a collection of targets $\Omega := \cup_{j=1} \mathcal{A}_j$. The shape of each target can be general and the boundary of each target can be combinations of absorbing and reflecting sections.  \label{Fig:schematic} }
\end{figure}

We consider the problem of describing the full arrival time distribution
of diffusing particles to complex absorbing sets in planar regions
$\RR^2\setminus\Omega$ (see Fig.~\ref{Fig:schematic}). For a particle
released from location $\bx_0\in\RR^2\setminus\Omega$, the central
quantities of interest are its occupation density $p(\bx,t;\bx_0)$ and
its survival probability
$P(t;\bx_0)=\int_{\RR^2\setminus\Omega}p(\bx,t;\bx_0)\,d\bx$ together
with the dependence on the number and location of absorbing sites. Our
new contributions are several new efficient numerical and asymptotic
methods to rapidly determine these quantities in the presence of complex
configurations of targets (cf.~Fig.~\ref{Fig:schematic}).

The general problem takes the form of a diffusion equation in $\bx \in
\RR^2\setminus\Omega$ where $\Omega := \cup_{j=1} \mathcal{A}_j$ is a collection of targets. We solve for the probability $p(\bx,t;\bx_0)$ of a particle originating at $\bx_0$ being free at position $\bx$ at time $t$. This probability density solves the exterior parabolic problem
\bsub\label{eqn:IntroP}
\begin{gather}
  \label{eqn:IntroP_a} 
  \frac{\partial p}{\partial t} = D \Delta p, \qquad
    \bx \in \RR^2\setminus\Om, \quad t>0; \qquad 
    p(\bx,0) = \delta (\bx-\bx_0), \quad \bx \in \RR^2\setminus\Om;\\[5pt]
  \label{eqn:IntroP_b} 
  p = 0, \quad \bx \in \partial\Omega_a; \qquad 
  D\nabla p \cdot \nn = 0, \quad \bx \in \partial\Omega_r ,
\end{gather}
\esub
where $D$ is the diffusivity of the particle and $\Omega$ is a
subset of $\RR^2$. The boundary $\pOm$ is partitioned into an
absorbing set $\Omega_a$ and its impermeable complement $\Omega_r$ where reflecting boundary conditions are applied. We choose $\nn$, the normal to the surface $\pOm$, to point
into the bulk. Some key quantities of interest for which rapid and
accurate determination is desirable include the fluxes over each target,
the splitting probabilities (likelihood particle encounters a certain
target first) and the arrival time distribution. The survival
probability
\begin{align*}
  P(t;\bx_0) = \int_{\RR^2\setminus\Omega} p(\bx,t;\bx_0) \, d\bx,
\end{align*}
is another important quantity obtained, together with its dependence on the number and location of targets.

First passage time problems and their variants appear in a variety of
disparate applications from cellular biology~\cite{Miles2020,Bressloff2021a,Bressloff2021b,BressNewby2013}, ecology~\cite{Kurella2015, Schenker2015,
Miller2017, Stepien2020}, and electrostatics~\cite{Chapman2015}. A few
comprehensive survey references can be found here~\cite{redner2001guide,
FPPA2014, berg, HolcmanReview2014}. Preceding works on the theory of
first arrival times to small absorbing sites have largely focussed on
determining moments such as the mean first passage time
(MFPT)~\cite{Pillay2010, ChevWard2010, HolcmanReview2014, LIW2021_a,
IW2021_b, IW2021_c, IW2021_d} and in some cases the
variance~\cite{Kurella2015,LWB2017,LTS2016}. In the present scenario of
an unbounded domain, these moments are not finite~\cite{redner2001guide,
FPPA2014} and we must therefore seek the full distribution of arrival
times. In the scenario where diffusion occurs in a bounded domain, the full arrival time distributions to small absorbing targets have been considered in two \cite{LTS2016,Bressloff2021a} and three \cite{Isaacson2013,Bressloff2021b} dimensions.  

In the planar (unbounded) scenario considered here, capture is guaranteed; however,
it may occur over very long timescales. This becomes apparent from the arrival time distribution $\chi_0(t)$ for a particle of diffusivity $D=1$, initially at distance $R$ from a target disc of unit radius centered at the origin. This distribution $\chi_0(t)$ and its large time behavior (see \cite{redner_2001} and Appendix \ref{sec:PlanarDist}) are given by
\begin{align}
  \label{eq:dist2Dreinsertion}
  \chi_0(t) = \frac{2}{\pi}\int_{0}^{\infty}\left[
    \frac{J_0(\omega)Y_0(\omega R)-J_0(\omega R)Y_0(\omega)}
    {Y_{0}^{2}(\omega)+J_{0}^{2}(\omega)}\right]\omega 
    e^{-\omega^2 t} \, d\omega = 
  \bigoh\left(\frac{1}{t\log^2 t}\right), \qquad t\to\infty,
\end{align}
where $J_0(z)$, $Y_0(z)$ are Bessel functions. The slow rate of decay in
the tail of this distribution reveals that very long arrival times are
typical (see Fig.~\ref{fig:2DarrivalCDF_a}). For example, when $R = 10$,
a particle with diffusivity $D=1$ still has an approximately $20\%$
chance of being free after $t = 10^{10}$. Equilibrium quantities (e.g.~splitting probabilities) therefore emerge on timescales that may not necessarily be the most biologically meaningful. For example, in applications such as a moth's search for a mate~\cite{Stepien2020},
or cellular signaling where a downstream event initializes as soon as a
molecule reaches a receptor~\cite{Miles2020}, the
statistics of particles which reach the target first are of most
interest. These \emph{extreme statistics} are governed by the behavior
of $P(t;\bx_0)$ for $t\ll1$ \cite{Lawley2020} and so it is necessary to
have methodologies for determining full distributions of arrival time
statistics.

In the present work, we outline several methods to solve
\eqref{eqn:IntroP}. First, in Sec.~\ref{sec:asy} we apply a Laplace
transform $\hat{p}(\bx,s)= \int_{t=0}^\infty p(\bx,t) e^{-st}\,dt$ to
\eqref{eqn:IntroP} to arrive at an elliptic problem of modified
Helmholtz type 
\bsub\label{eq:LT2D}
\begin{gather}
  \label{eq:LT2D_b} 
  D\Delta\hat{p} - s\, \hat{p} = -\delta(\bx-\bx_0),\qquad
    \bx\in\RR^2\setminus\Omega;\\[5pt]
  \label{eq:LT2D_a}
  \hat{p} = 0 \quad \bx \in \Omega_a; \qquad 
    \nabla \hat{p} \cdot \textbf{n} = 0 \quad \bx \in \Omega_r.
\end{gather}
\esub
In the limit of well-separated traps, we solve the resulting transform
problem in terms of an asymptotic expansion where solutions are obtained
in terms of modified Helmholtz Green's function. This methodology was
originally developed in~\cite{LTS2016} and recently applied
in~\cite{Bressloff2021a,Bressloff2021b} to the study of first passage
times of particles with resetting. The Laplace transform is inverted
numerically by Talbot quadrature~\cite{Abate2006} resulting in a hybrid
numerical-asymptotic method. In Sec.~\ref{sec:BIE}, we take a similar
approach, but replace the asymptotic solution of~\eqref{eq:LT2D} with a
layer potential representation. This results in a boundary integral
equation that is solved numerically using a collocation method.

In Sec.~\ref{sec:KMC} we develop a particle based kinetic Monte-Carlo
(KMC) method that evaluates solutions of \eqref{eqn:IntroP} by dividing
the sojourn of particles into projection steps where exact solutions are
available \cite{LBS2018,OKMascagni2004}. This offers a rapid, accurate
and easy to implement method for the solution of \eqref{eqn:IntroP} in very
general geometries. In Sec.~\ref{sec:results} we demonstrate the
applicability of these methods on a variety of examples. In particular,
we provide numerical validation of previously derived homogenization
theories and find them to be highly effective in reproducing the arrival
time distributions. We also investigate the time dependent fluxes into the
targets which are observed to converge very slowly to the static splitting 
probabilities that describe the relative flux into each target. This suggests that
a relevant physical or biological time scale should be considered before using
 receptor arrival information to make inferences on environmental conditions.

% 2D problem: analysis
\section{Asymptotic description of arrival times of particles diffusing
in $\RR^2$}
\label{sec:asy}
In this section we use matched asymptotic expansions to derive an
approximation for the density of a particle diffusing in $\RR^2\setminus\Omega$
in the presence of well-separated target sites. The assumption is that
there are $N$ targets $\mathcal{A}_j$ with centers $\{\bx_j\}_{j=1}^N$ so that
the collection of target sites are described by
\begin{align}\label{eq:2DAbsSet}
\Omega = \bigcup_{j=1}^N  \mathcal{A}_j \left (\frac{\bx -\bx_j} {\eps} \right ), \qquad \bx = (x,y),
\end{align}
where $\eps$ is a parameter controlling the extent of the targets and enforces the well-separated condition as $\eps\to0$.
The geometry of individual targets can be quite general.

The aim is to solve the solution of the initial-boundary value problem
\eqref{eqn:IntroP} and determine the free probability
$P(t)=\int_{\Omega} p(\bx,t)\,d\bx$ together with the capture time
density $C(t)= -P'(t)$. The first step \cite{LTS2016} in the analysis is to define, for
$s\in\CC$, the Laplace transform $\hat{p}(\bx,s)= \int_{t=0}^\infty
p(\bx,t) e^{-st}\,dt$ that solves \eqref{eq:LT2D}. The mixed boundary
conditions \eqref{eq:LT2D_a} indicate that the target boundary may have
a combination of absorbing or reflecting components so that
$\partial\Omega = \Omega_a\cup\Omega_r$. In the absence of the target
set $\Omega$, the solution $\hat{p}$ of \eqref{eq:LT2D} is defined in
terms of the free space modified Helmholtz Green's function
$G_h(\bx;\bxi,s)$
\bsub\label{eqG}
\begin{gather}
\label{eqG_a} D\Delta G_h-s\, G_h=-\delta(\bx-\bxi), \qquad
  \bx\in\RR^2\setminus\{\bxi\};\\[5pt]
\label{eqG_b} G_h(\bx;\bxi,s)=\frac{1}{2\pi D} K_0(\sqrt{s/D}|\bx-\bxi|),
  \qquad \int_{\RR^2} G_h(\bx;\bxi,s)\,d\bx=\frac{1}{s}; \\[5pt]
\label{eqG_c} G_h(\bx;\bxi,s) \sim -\frac{1}{2\pi D} \log |\bx-\bxi| + R_h (s) + \littleoh(1), \qquad \bx\to\bxi.
\end{gather}
\esub
Here $R_h(s)$ is the regular part of $G_h(\bx;\bxi,s)$ at the source. The small argument asymptotics $K_0(z)\sim-\log(z)+\log2-\gamma_e$ as $z\to0$ give this self-interaction term to be
\begin{equation}
 R_h(s)=\frac{1}{2\pi D}\left(\log2-\gamma_e-\log\sqrt{s/D}\right),
\end{equation}
where $\gamma_e\approx 0.5772$ is the Euler-Mascheroni constant.

In the limit of well-separated absorbers $\eps\to0$, we employ a matched
asymptotic analysis to replace each target \eqref{eq:LT2D_a} by
effective singularity conditions. To establish this singularity
condition, the following change of variables is introduced near the
$j^{\mathrm{th}}$ absorber
\begin{equation}
 \by=\frac{\bx-\bx_j}{\eps}, \qquad v(\by) = \hat{p}(\bx_j + \eps \by).
\end{equation}
In these coordinates, the transformed equation \eqref{eq:LT2D_b} is
$D\Delta_{\by}v - s \eps^2 v= - \eps^2 \delta(\bx-\bx_0)$. In addition
to the limit $\eps\to0$, we additionally consider the case $s\eps^2
\ll1$ which is valid provided $s$ is not too large. The limit
$s\to\infty$ corresponds to $t\to0$ and therefore we cannot expect good
agreement for arbitrarily short times. With these points in mind, we
continue by considering the local solution $v(\by) = v_j (\by)+
\bigoh(\eps)$ near the target where $v_j(\by)$ satisfies the exterior
problem
\bsub\label{eqn:LogCap}
\begin{gather}
\Delta_{\by}v_j=0,\qquad \by\in\RR^2/\A_j,\\[5pt]
v_j(\by)\sim \log|\by|-\log d_j + \frac{ {\textbf p} \cdot \by }{|\by|^2} + \cdots ,\quad |\by|\to\infty.
\end{gather}
\esub
Here the parameter $d_j$ is the logarithmic capacitance and depends on
the shape of $\A_j$ and the boundary conditions applied to it. In
section \ref{sec:LogCap} we give an overview of many scenarios in which
$d_j$ can be calculated. The behavior of $v_j(\by)$ at infinity gives
the matching condition for $\hat{p}(\bx)$ as $\bx$ approaches $\bx_j$.
That is,
\begin{align}\label{eq:guage}
\hat{p}(\bx,s)\sim S_j \nu_j
  \left(\log\left|\frac{\bx-\bx_j}{\eps}\right|-\log d_j\right)=S_j
  \nu_j \log|\bx-\bx_j|+S_j,\quad \bx\to\bx_j; \qquad \nu_j = -\frac{1}{\log \eps d_j},
\end{align}
where $S_j$ is a strength term to be determined in terms of $(N+1)$ parameters $(s,\bs{\nu})=(s,\nu_1,\ldots,\nu_N)$. Therefore, in the outer region away from targets, we pose the asymptotic expansion
\begin{align*}
  \hat{p}(\bx;s) = \hat{p}_0(\bx;s,\bs{\nu}) + \eps \,\hat{p}_1(\bx;s,\bs{\nu}) + \littleoh(\eps), \qquad \eps\to0.
\end{align*}
The leading order solution $\hat{p}_0(\bx;s,\bs{\nu})$ \lq\lq sums-the-logs\rq\rq\ and is accurate to all logarithmic orders. The correction term $\hat{p}_1$, which we do not explicitly determine, describes how target orientation influences capture and can be found following methods outlined in \cite{LTK2015}. The leading order problem satisfies 
\bsub\label{equ2}
\begin{gather}
D\Delta\hat{p}_0 -s\, \hat{p}_0=-\delta(\bx-\bx_0),\quad \bx\in\RR^2/\{\bx_1,\ldots,\bx_N\},\\[5pt]
\hat{p}_0\sim S_j \nu_j \log|\bx-\bx_j|+S_j,\quad \bx\to\bx_j, \qquad j = 1,\ldots, N.
\end{gather}
\esub
The solution $\hat{p}_0(\bx, s)$ of \eqref{equ2} is described in terms of the modified Helmholtz Green's function \eqref{eqG} as
\begin{equation}
\hat{p}_0(\bx,s)= G_h(\bx;\bx_0,s)- 2\pi D\sum_{j=1}^N S_j \nu_j G_h(\bx;\bx_j,s).
\label{usol}
\end{equation}
The coefficients $S_j$ can be determined by equating the regular parts in (\ref{equ2}b) and in (\ref{usol}). That is,
\begin{equation}
S_j= 
  G_h(\bx_j;\bx_0,s)- 2\pi D\left[S_j\nu_jR_h(s)+\sum_{\substack{i=1\\
  i\neq j}}^NS_i\nu_iG_h(\bx_j;\bx_i,s)\right], \qquad j = 1,\ldots,N.
\label{eqA}
\end{equation}
\bsub\label{sysmain}
In summary, we have that the transform equation \eqref{eq:LT2D} has asymptotic solution $p(\bx;s) \sim p_0(\bx;s)+\cdots$ as $\eps\to0$ where
$p_0(\bx;s)$ satisfies
\begin{equation}
\hat{p}_0(\bx,s)= G_h(\bx;\bx_0,s)- 2\pi D\sum_{j=1}^N S_j \nu_j G_h(\bx;\bx_j,s).
\end{equation}
The strengths $\{S_j\}_{j=1}^N$ satisfy \eqref{eqA} which can be represented in compact matrix form as
\begin{equation}
 (\mathcal{I} + 2\pi D\,  \mathcal{G}_h \mathcal{V})\, \textbf{S}=\bg_0, \qquad [\mathcal{G}_h]_{i,j}= \left\{ \begin{array}{rl} R_h(s),\quad i= j, \\[5pt]  G_h(\bx_i;\bx_j,s),\quad i\neq j. \end{array} \right. \qquad [\mathcal{V}]_{i,j}= \left\{ \begin{array}{rl} \nu_i \quad i= j, \\[5pt]  0,\quad i\neq j, \end{array} \right. 
\label{matrixG}
\end{equation}
where $\mathcal{I}\in\RR^{N\times N}$ is the identity,
$\textbf{S}\in\RR^N$ and $ \bg_0\in\RR^N$ are given by
\begin{align}
\textbf{S}= [S_1,S_2,\ldots,S_N]^T,\qquad
 \bg_0=[G_h(\bx_1;\bx_0,s), G_h(\bx_2;\bx_0,s), \ldots,G_h(\bx_N;\bx_0,s)]^T.
\end{align}
\esub
The matrix $\mathcal{G}$ describes the interactions between targets and
their competition for flux while the vector $\bg_0$ reflects the
influence of the initial location on each of the targets.  The vector
$\widehat{\mathcal{J}} = 2\pi D \mathcal{V}\textbf{S}(s)$ describes the transformed fluxes through each of the targets and is obtained from solving the linear system \eqref{sysmain}. 

At this stage we calculate additional quantities of interest, namely the survival probability $P(t)$ and arrival time distribution $C(t)$. Using equation \eqref{usol} and \eqref{eqG_c}, the Laplace transform $\hat{P}(s)$ of the free probability is given by
\bsub\label{2DResults}
\begin{equation}\label{2DResults_P}
 \hat{P}(s)=\int_{\Omega}
  \hat{p}_0(\bx,s)\,d\bx= \int_{\Omega} G_h(\bx;\bx_0,s)d\bx - 2\pi D\sum_{j=1}^N S_j \nu_j \int_{\Omega} G_h(\bx;\bx_j,s)d\bx =\frac{1}{s}\left[1-2\pi D\sum_{j=1}^N \nu_j S_j(s)\right].
\end{equation}
The relationship $C(t)=-P'(t)$, yields that the Laplace transform of the arrival time distribution $C(t)$ is
\begin{equation}\label{2DResults_C}
\hat{C}(s)=-[s\hat{P}(s)-P(0)]=-s\hat{P}(s)+1=2\pi D\sum_{j=1}^N \nu_j S_j(s).
\end{equation}
\esub
\begin{figure}[htbp]
\centering
\includegraphics[width = 0.5\textwidth]{TalbotSchematic.pdf}
\caption{A schematic of the Talbot curve \eqref{eq:Talbot} for parameter values $\sigma=0$, $\mu = 2$, $\nu = 0.5$. The red line indicates the singularities along the negative real line arising from the $\sqrt{s}$ singularity. \label{fig:Talbot} }
\end{figure}

\subsection{Inverse Laplace Transform}\label{sec:invLT}

To obtain $P(t)$ and $C(t)$ defined by \eqref{2DResults}, the inverse Laplace transform
\begin{equation}\label{eqn:numInvLT}
P(t) = \frac{1}{2\pi i} \int_{\Gamma_B} e^{st} \hat{P}(s)  \, \textrm{d}s,
\end{equation}
must be evaluated where $\Gamma_B$ is the Bromwich contour $\Gamma_B = \{ \gamma + i y \, | \, -\infty<y<\infty \}$. The parameter $\gamma$ is chosen so that all singularities of $\hat{P}(s)$ lie to the left of Re($s)=\gamma$. In the present scenario associated with diffusive motion, the singularities of $\hat{P}(s)$ lie along the negative real axis due to the branch cut of $\sqrt{s}$. Rapid and effective numerical evaluation of \eqref{eqn:numInvLT} can be achieved by deforming the contour around Re($s)=0$ since the integrand of \eqref{eqn:numInvLT} decays very rapidly for Re($s)<0$. The Talbot contour $\Gamma_T$ is a family of deformations (see Fig.~\ref{fig:Talbot})  to $\Gamma_B$ where
\begin{equation}\label{eq:Talbot}
\Gamma_T = \{ \sigma + \mu (\theta \cot \theta + \nu i \, \theta)  \, | \, -\pi<\theta<\pi \},
\end{equation}
and $\sigma, \mu$ and $\nu$ are parameters that control the curve shape
\cite{Weideman2006,Trefethen2006}. Rapid and accurate evaluation of the
inverse Laplace transform is then achieved by applying the midpoint rule
on this curve.

\subsection{Logarithmic capacitance for various shapes}\label{sec:LogCap}

The asymptotic solution \eqref{2DResults} encodes the geometry of each target $\mathcal{A}_j$ into the logarithmic capacitance $d_j$, determined by the solution of \eqref{eqn:LogCap}. Here we discuss the determination of $d_j$ and briefly recap known results for regular geometries and simple boundary conditions.

\paragraph{Regular geometric shapes} For simple shapes such as circles, ellipses, triangles and squares with all absorbing perimeters, the logarithmic capacitance is known exactly. A list of these quantities, reproduced from \cite{Kurella2015}, is included in Table~\ref{tab:LogCap}.

\begin{table}[H]
    \centering
    \begin{tabular}{c|c}
    \toprule
         {Shape of $\A_j$} & Logarithmic capacitance $d_j$  \\
         \midrule
         \midrule
         circle of radius $a$ & $d_j = a$\\[5pt]
         ellipse, semi-axes $a,b$ & $d_j = \frac{a+b}{2}$\\[5pt]
         equilateral triangle, side-length $h$ & $d_j = \frac{\sqrt{3}\Gamma(\frac13)^3 h}{8\pi^2} \approx 0.422 h$\\[5pt]
         isosceles right triangle, side-length $h$ & $d_j = \frac{3^{3/4}\Gamma(1/4)^2 h}{2^{7/2}\pi^{3/2}} \approx 0.476 h$\\[5pt]
         square, side-length $h$ & $d_j = \frac{\Gamma(1/2)^2 h}{4\pi^{3/2}} \approx 0.590 h$\\
         \bottomrule
    \end{tabular}
    \bigskip
    \caption{The logarithmic capacitance of some simple geometries with
    absorbing boundary conditions, reproduced from \cite{Kurella2015}.
    \label{tab:LogCap}}
\end{table}

\paragraph{Partially absorbing disk: Single window} For a circular trap that is absorbing except for the reflecting portion $\theta \in(-\sigma,\sigma)$, the problem \eqref{eqn:LogCap} may be expressed in polar coordinates $(r,\theta)$ as
\bsub\label{eqn:vc_partial}
\begin{gather}
    \Delta v = 0, \quad r \geq 1, \quad \theta\in(-\pi,\pi);\\[5pt]
    \partial_r v = 0, \qquad r = 1, \quad \theta \in(-\sigma,\sigma); \qquad 
    v = 0,\qquad r = 1, \quad \theta \in(\sigma,\pi)\cup (-\pi,-\sigma);\\[5pt]
    v \sim \log r - \log d_c + \bigoh\big(r^{-1}\big), \quad r \to \infty.
\end{gather}
\esub
The separable solution of \eqref{eqn:vc_partial} takes the form of the cosine series 
\begin{align*}
  v(r,\theta) = \log r + \frac{a_0}{2} + 
    \sum_{n=1}^{\infty}\frac{a_n}{r^n}\cos n\theta,
\end{align*}
where the coefficients $\{a_n\}_{n=0}^{\infty}$ satisfy the dual
trigonometric series
\bsub\label{eqn:dualTrig}
\begin{align}
    \frac{a_0}{2} + \sum_{n=1}^{\infty} a_n \cos n \theta &= 0, \quad \theta\in(\sigma,\pi),\\
    \sum_{n=1}^{\infty} n a_n \cos n \theta &= 1, \quad \theta\in(0,\sigma).
\end{align}
\esub
The solution of \eqref{eqn:dualTrig} was determined in \cite{LTK2015} from an integral equation theory which reveals that
\begin{equation}\label{eqn:a0split}
   -\log d_c = \frac{a_0}{2} = \frac{\sqrt{2}}{\pi} \int_0^{\sigma} \frac{u \sin \frac{u}{2}}{\sqrt{\cos u - \cos \sigma}}\, \mathrm{d}u.
\end{equation}
For the half absorbing case ($\sigma = \pi/2$), $a_0 = 2\log2$, while in the singular mostly reflecting limit $\sigma\ll1$, it can be determined that $a_0 \sim -4\log\frac{\sigma}{2}$. In other scenarios, the integral \eqref{eqn:a0split} is readily evaluated by quadrature.

\paragraph{Partially absorbing disk: Multiple windows}
In the scenario of a circular trap with $N$ small absorbing windows of length $\sigma$ centered at points $\{\by_j \}_{j=1}^N$, it was shown in \cite{LTK2015} that as $\sigma\to0$
\bsub\label{eq:multi_windows}
\begin{equation}\label{eq:multi_windows_a}
-\log d_c = \frac{a_0}{2} \sim -\frac{2}{N} \log\frac{\sigma}{4} - \frac{2}{N^2}\sum_{i=1}^N \log \prod_{\substack{j=1 \\ j\neq i}}^N |\by_i - \by_j|.
\end{equation}
For windows centered at roots of unity $\by_j = \left(\cos \frac{2\pi
j}{N},\sin \frac{2\pi j}{N}\right)$, and for $N\sigma <1$,
\eqref{eq:multi_windows_a} reduces to
\begin{equation}
  \label{eq:multi_windows_b}
  -\log d_c = \frac{a_0}{2} = -\frac{2}{N} \log\frac{\sigma}{4} - 
    \frac{2}{N^2}(N\log N)= -\frac{2}{N} \log \frac{\sigma N}{4}.
\end{equation}
\esub
\paragraph{Partially absorbing disk: Homogenization limit}
The results \eqref{eq:multi_windows} can be used to identify an
homogenization limit as $N\to \infty$ and $\sigma\to0$. The  absorbing
fraction $f$ is defined through $N\sigma = 2\pi f$ and the homogenized
logarithmic capacitance problem satisfies
\bsub\label{eqn:homogenization}
\begin{gather}
\label{eqn:homogenization_a} \Delta v_h=0,\quad r>1; \qquad \sigma \frac{\partial v_h}{\partial r} + \kappa(f) v_h = 0, \quad r = 1;\\[5pt]
\label{eqn:homogenization_b} v_h\sim \log r-\log d_c + \cdots ,\quad r\to\infty.
\end{gather}
In the dilute limit $f\ll1$, the homogenized parameters were identified in \cite{LTK2015} to be
\begin{equation}\label{eqn:homogenization_c}
  \kappa(f) = -\frac{\pi f}{\log \frac{\pi f}{2}}, \qquad  \log d_c = \frac{\sigma}{\kappa(f)}\, .
\end{equation}
\esub
In Section~\ref{sec:resultsHomogenization}, we show numerical results
that validate this homogenized formulation and demonstrate that it is highly accurate in predicting the arrival time statistics of the full problem.

\paragraph{The logarithmic capacitance for a two trap cluster}
In the case of two circular traps separated by distance $\ell$, it was
derived in \cite{Kurella2015} (see also
\cite{LIW2021_a,IW2021_b,Strieder08,Strieder12}) from an expansion in
bi-polar coordinates that 
\begin{equation}\label{eq:TwoTrapCluster}
\log d_c = \frac{1}{2}\log(\ell^2 - 4) - \frac{\beta}{2} + \sum_{k=1}^{\infty}\frac{e^{-k\beta}}{k\cosh(k\beta)}, \qquad \beta = \cosh^{-1} \left(\frac{\ell}{2}\right).
\end{equation}

\paragraph{Numerical evaluation of the logarithmic capacitance for general configurations}

For very general configurations of clusters, the logarithmic capacitance problem \eqref{eqn:LogCap} can be obtained numerically by a boundary integral approach \cite{Dijkstra2008}. Another approach developed in \cite{LIW2021_a,IW2021_b}, and based on works \cite{Chapman2015,trefethen_2018}, is to develop a series solution of \eqref{eqn:LogCap} followed by a least squared method to obtain the unknown coefficients. For the case of $m$ circular absorbers with centers $\{c_j\}_{j=1}^m$, the series takes the form
\begin{equation}\label{eqn:complexSeries}
    v(z) = - \log d_c + \sum_{j=1}^m e_j \log|z - c_j| + \sum_{j=1}^m
    \sum_{k = 1}^n \left( a_{j k} \text{Re}(z-c_j)^{-k} + b_{j k}
    \text{Im} ( z- c_j) ^{-k} \right), \quad z\in\CC.
\end{equation}
The constants $d_c, e_j,a_{jk}, b_{jk}$ are to be determined while $\sum_{j=1}^{m}e_j=1$ enforces the far field behavior $v \sim \log |z| $ as $|z|\to\infty$. A system for the unknown constants is formed by evaluating \eqref{eqn:complexSeries} at a collection of boundary points along which $v = 0$. A similar methodology was used in \cite{LTK2015} to solve a truncated version of the dual trigonometric series \eqref{eqn:dualTrig} by evaluation at a set of boundary values followed by least squared solution. 

% inverse Laplace transform

\section{Boundary integral equation description of arrival times of
particles diffusing in $\RR^2$}
\label{sec:BIE}
An alternative approach to matched asymptotics for
solving~\eqref{eqn:IntroP} is to use an integral equation approach.
Integral equations are a natural choice for unbounded complex domains
such as the one in Fig.~\ref{Fig:schematic} since they easily resolve
complex geometries while automatically satisfying the far-field boundary
conditions. Others have applied integral equation methods to
solve~\eqref{eqn:IntroP} using the full space-time heat
kernel~\cite{gre-li2009} or by discretizing in time and solving the
resulting elliptic PDE with an integral equation
formulation~\cite{kro-qua2010, cha-kre1997}. However, these approaches
have several challenges that include maintaining long time
histories and computing volume integrals. We take a new approach by
solving for the Laplace transformed variable $\hat{p}(\xx,s)$ that
satisfies~\eqref{eq:LT2D}. We only consider the case where $\bd\Omega$
is absorbing so that the boundary condition is Dirichlet and
homogeneous.

We begin by writing $\hat{p}$ as the sum of a particular and homogeneous
solution of~\eqref{eq:LT2D}
\begin{align*}
  \hat{p}(\xx,s) = G_h(\xx;\xx_0,s) + \hat{p}^{H}(\xx,s),
\end{align*}
where $G_h$ is the free space modified Helmholtz Green's
function~\eqref{eqG}. Using the boundary condition~\eqref{eq:LT2D_a},
$\hat{p}^H$ satisfies the homogeneous PDE
\begin{subequations}
\label{eqn:homoPDE_BC}
\begin{alignat}{3}
  D \Delta \hat{p}^{H} - s \hat{p}^H &= 0, &&\xx \in \Omega, 
  \label{eqn:homoPDE} \\[5pt]
  \hat{p}^{H}(\xx) &= f(\xx), \quad && \xx \in \bd\Omega,
  \label{eqn:homoBC}
\end{alignat}
\end{subequations}
where $f(\xx) = -G_h(\xx;\xx_0,s)$. We represent the solution
of~\eqref{eqn:homoPDE_BC} with the double-layer potential
\begin{align*}
  \hat{p}^{H}(\xx) = \DD[\sigma](\xx) = 
    \int_{\bd\Omega} \pderiv{}{{\nn}_\yy}
    G_{h}(\xx;\yy,s) \sigma(\yy) \, ds_\yy, \quad \xx \in \Omega,
\end{align*}
where $\sigma$ is an unknown density function. We remind the reader that
the unit normal $\nn_\yy$ points into the bulk. To satisfy the boundary
condition~\eqref{eqn:homoBC}, the density function must solve
\begin{align}
  \label{eqn:SKIE}
  f(\xx) = \frac{1}{2} \sigma(\xx) + \DD[\sigma](\xx), 
    \quad \xx \in \bd\Omega.
\end{align}
A numerical solution of the second-kind integral
equation~\eqref{eqn:SKIE} is formed by discretizing $\bd\Omega$ at $N$
quadrature points and approximating the integrals with the trapezoid
rule. The resulting linear system is
\begin{align}
  \label{eqn:discretizedSKIE}
  f(\xx_i) = \frac{1}{2} \sigma(\xx_i) + 
    \frac{1}{2\pi} \sum_{j=1}^{N} G_{h}(\xx_i;\xx_j,s) 
    \sigma(\xx_j) \, \Delta s_j, \quad i=1,\ldots,N,
\end{align}
and the diagonal term of this linear system is replaced with the
limiting value
\begin{align*}
  \lim_{\substack{\yy \rightarrow \xx \\ \xx \in \bd\Omega}}
    G_h(\xx;\yy,s) = -\frac{1}{2}\kappa(\xx),
\end{align*}
where $\kappa(\xx)$ is the curvature of $\bd\Omega$. The linear
system~\eqref{eqn:discretizedSKIE} is solved with the generalized
minimal residual method (GMRES), and since it is the discretization of a
second-kind integral equation, the number of required iterations is
mesh-independent. This method to solve for $\hat{p}(\xx,s)$ is coupled
with the inverse Laplace transform~\eqref{eqn:numInvLT}, where we use
the same Talbot contour illustrated in Fig.~\ref{fig:Talbot}.

The flux at point $\xx \in \bd\Omega$ and time $t$ is $\JJ(\xx,t) =
\pderiv{}{\nn}p(\xx,t)$. Since we write $\hat{p}(\xx,s)$ as the sum of a
fundamental solution and a homogeneous solution, we compute the flux of
these terms individually, and the flux due to the fundamental solution
is computed analytically. The flux due to the homogeneous solution
$\hat{p}^H$ is
\begin{align}
  \label{eqn:DLPflux}
  \widehat{\JJ}^H(\xx;s) = \pderiv{}{\nn_\xx} \DD[\sigma](\xx) 
  = \int_{\bd\Omega} \pderiv{}{\nn_\xx}\pderiv{}{\nn_\yy}
  G_h(\xx;\yy,s) \sigma(\yy) \, ds_\yy,
\end{align}
which needs to be estimated with quadrature. The trapezoid rule or
odd-even integration~\cite{sid-isr1988} cannot be used to
approximate~\eqref{eqn:DLPflux} because
\begin{align*}
  \pderiv{}{\nn_\xx}\pderiv{}{\nn_\yy} G_h(\xx;\yy,s) \sim 
    \|\xx - \yy\|^{-2}.
\end{align*}
To formulate the normal derivative of the double-layer potential with a
tractable integrand, we first add and subtract the leading order
asymptotics of $G_h$ described in~\eqref{eqG_c}. That is,
$\widehat{\JJ}^H(\xx;s) = I_1 - I_2$ where
\begin{align*}
  I_1 &= \int_{\bd\Omega} 
    \pderiv{}{\nn_\xx}\pderiv{}{\nn_\yy}
    \left(G_h(\xx;\yy,s) + \frac{1}{2\pi} \log|\xx - \yy|\right) 
    \sigma(\yy) \, ds_\yy,  \\
  I_2&=\int_{\bd\Omega} \pderiv{}{\nn_\xx}\pderiv{}{\nn_\yy}
    \left(\frac{1}{2\pi} \log|\xx - \yy|\right) \sigma(\yy) \, ds_\yy.
\end{align*}
The singularity of the integrand in $I_1$ behaves as $\|\xx -
\yy\|^{-1}$, and odd-even integration can be applied. The integral
$I_2$ is further decomposed as
\begin{align*}
  I_2 &= \frac{1}{2\pi}\int_{\bd\Omega} 
    \pderiv{}{\nn_\xx}\pderiv{}{\nn_\yy}\log|\xx - \yy|
    (\sigma(\yy) - \sigma(\xx)) \, ds_\yy 
    - \frac{\sigma(\xx)}{2\pi} \int_{\bd\Omega} 
      \pderiv{}{\nn_\xx}\pderiv{}{\nn_\yy} \log|\xx - \yy| \, ds_\yy.
\end{align*}
The second integral in this expression is the normal derivative of a
constant function, and therefore is zero. The remaining integral has an
integrand with a singularity that also behaves as $\|\xx - \yy\|^{-1}$,
and odd-even integration can be applied.

Having developed a quadrature method to compute $\widehat{J}^H(\xx;s)$,
the point-wise flux can be computed at time $t$ by applying the midpoint
rule along the Talbot contour in Fig.~\ref{eq:Talbot}. Then, the total
flux into $\bd\Omega$ can easily be computed by applying the trapezoid
rule to
\begin{align*}
  S(t) = \int_{\bd\Omega} \JJ(s,t) \, ds.
\end{align*}

\section{Particle based Kinetic Monte Carlo simulations}\label{sec:KMC}

Monte Carlo simulations provide a valuable tool for numerically estimating the distribution of capture times of diffusing particles for problems such as \eqref{eqn:IntroP}  and have been used extensively \cite{berez2014, Northrup1988,Litwin1980,Opplestrup2006}. Trajectories associated to the density \eqref{eqn:IntroP} can be constructed through the discretization
\begin{equation}\label{iter1}
\bx(t+ dt) = \bx(t) + \sqrt{2D dt} \, Z, \qquad \bx(0) = \bx_0,
\end{equation}
where $Z\sim N(0,1)$. The sequence of small displacements \eqref{iter1} terminates when the particle encounters the absorbing surface $\partial\Omega_a$. The algorithm is repeated for many particles (millions or even billions) to sample the capture time distribution. This approach is flexible and easy to implement but hampered by a set of problems. 

If a fixed stepsize $dt$ is adopted, errors are introduced at that lengthscale that accrue near boundaries. First, for a capture event in the interval $(t,t+dt)$, we typically choose $t+dt$ as the arrival time which is an overestimate. Second, trajectories drawn from \eqref{iter1} will necessarily miss some encounters with boundaries and therefore overestimate the hitting time. Another challenge is that capture problems associated with \eqref{eqn:IntroP} are notorious for their fat-tailed distributions, i.e., a significant fraction of realizations undergo long excursions before capture. A key component of any
efficient method is adaptivity in stepsize since a trajectory of \eqref{iter1} simulated with a fixed step method will take a very long time to reach an absorbing site.

\subsection{Kinetic Monte Carlo (KMC) method for simulation of planar diffusion to absorbers}

Decreasing the step size in an adaptive manner based on distance to
target can ameliorate these issues. The Kinetic Monte Carlo (KMC) method
\cite{BatBerShv2003} maximizes this opportunity by advancing the
diffusion process in a spatial stepsize corresponding to the distance to
the target, $d(\bx_0,\partial\Omega)$. The geometry of each step can take
many forms, but it should be chosen such that the details of the sojourn
can be rapidly and accurately sampled from closed form expressions.
Similar ideas have been employed in $N$-body simulations of kinetic
gases \cite{Opplestrup2006} and chemical reactions \cite{Wu2006}. In
this paper we describe implementation details and rudimentary analysis
of such a scheme that can handle complex geometries and mixed boundary
conditions. This method completely bypasses the need to advance
particles based on discretized steps such as~\eqref{iter1}.

\begin{figure}[htbp]
\centering
\subfigure[Stage I: walk on sphere step.]{\includegraphics[width=0.3\textwidth]{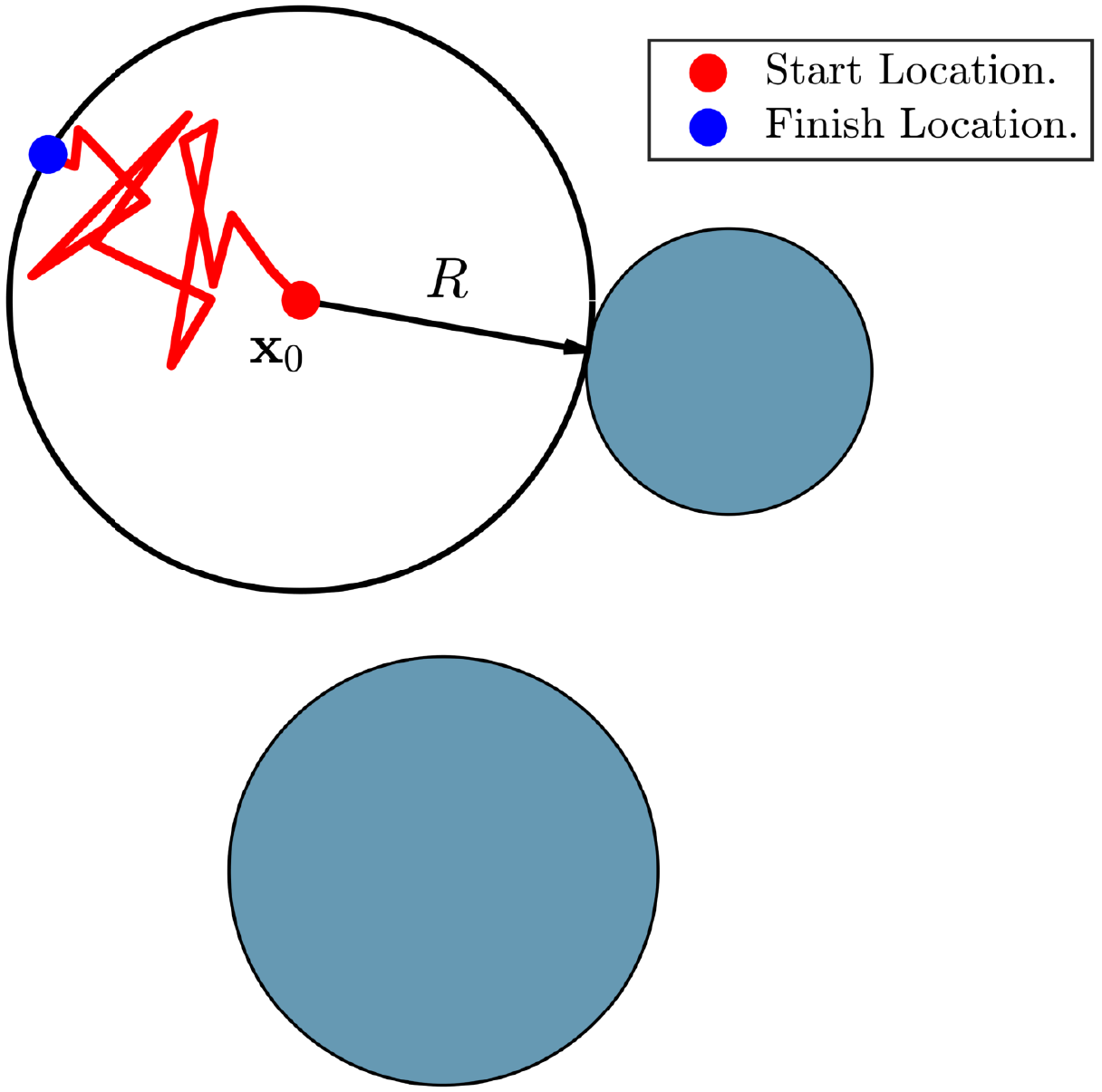}\label{fig:WoS1a}}\qquad
\subfigure[Stage II: reinsertion step.]{\includegraphics[width=0.4\textwidth]{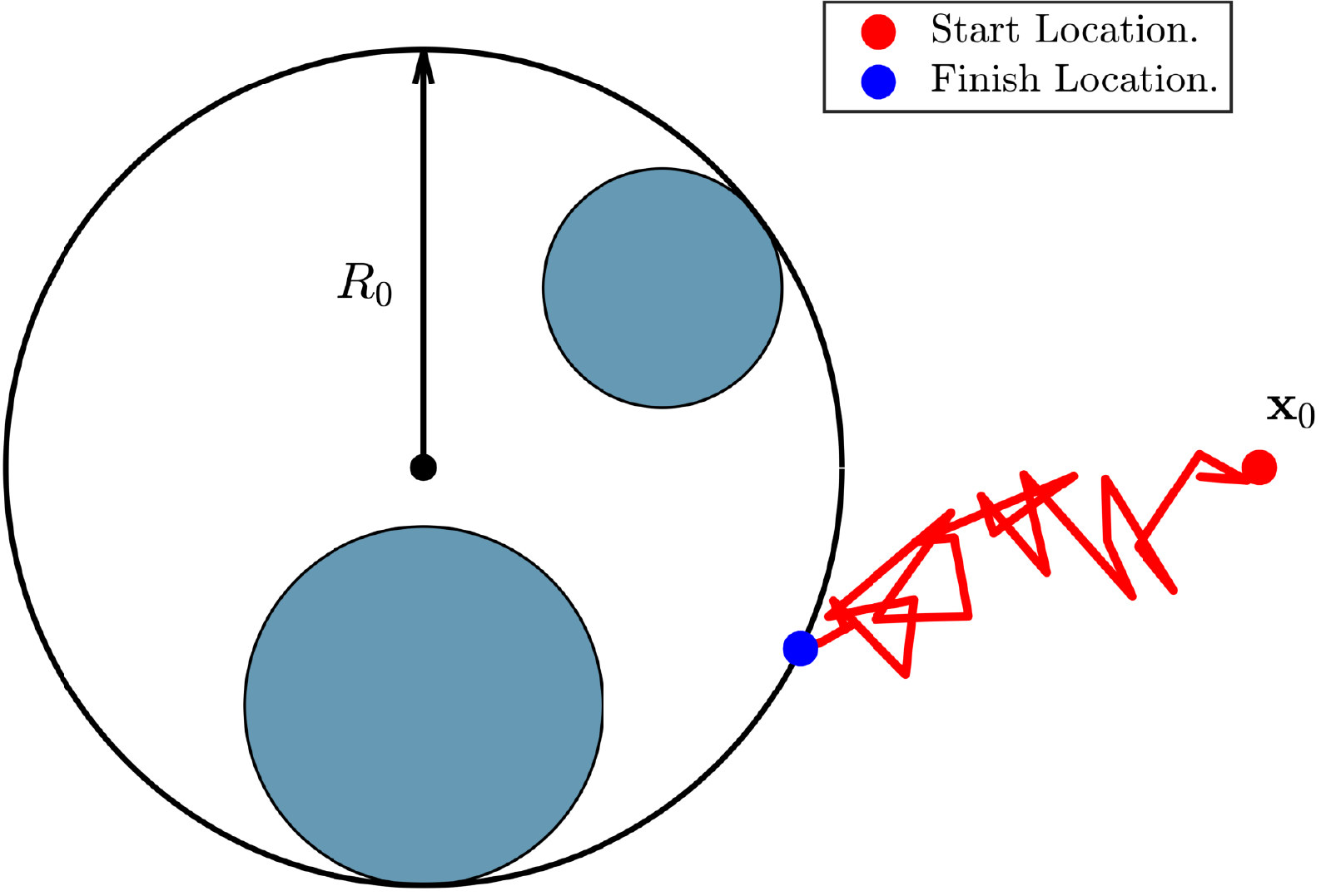}\label{fig:WoS1b}}
\subfigure[Stage III: square projection step.]{\includegraphics[width=0.4\textwidth]{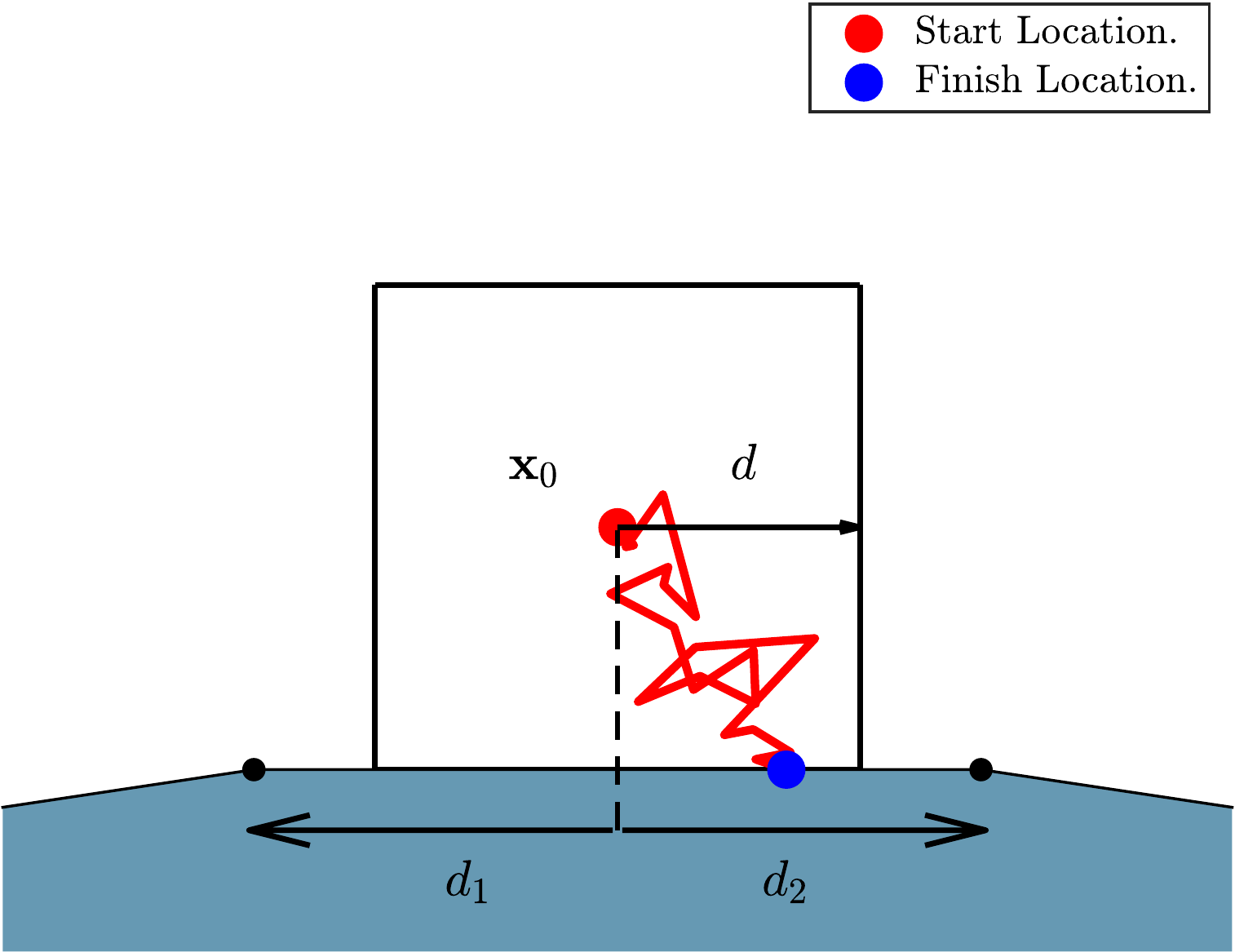}\label{fig:WoS1c}}\qquad
\subfigure[Stage IV: reflection step.]{\includegraphics[width=0.4\textwidth]{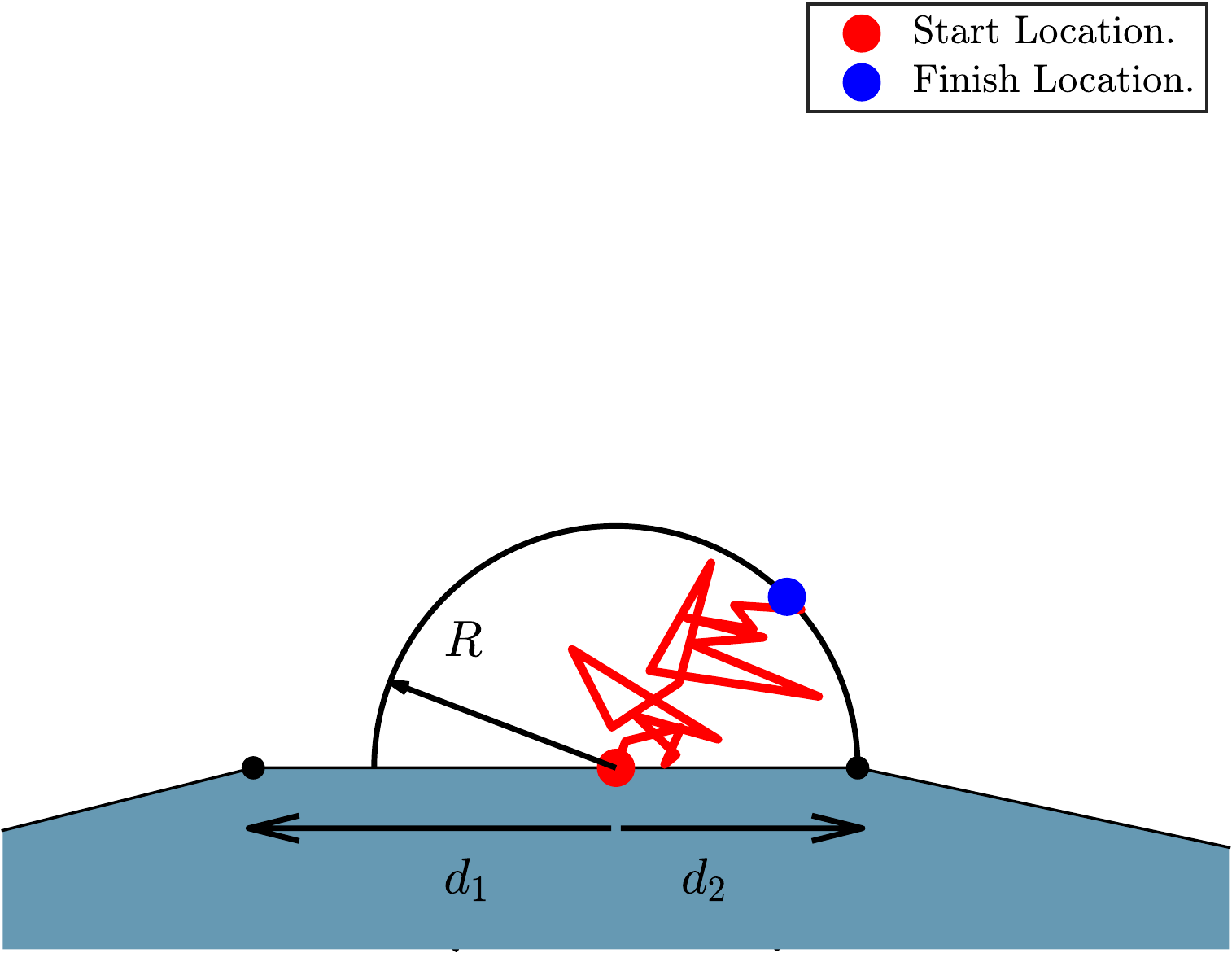}\label{fig:WoS1d}}
\caption{Schematic of the four stages of the planar KMC algorithm. \label{fig:WoS1}}
\end{figure}

\paragraph{Setup} We adopt a piecewise linear representation for the
boundary $\partial\Omega$ of the target set $\Omega$ based on vertexes
with $M$ straight edges $\partial\Omega_j$ so that $\partial \Omega =
\bigcup_{j=1}^M\partial\Omega_j$. On each boundary edge
$\partial\Omega_j$, we precalculate midpoints, unit normal vectors and
associate either Neumann or Dirichlet boundary conditions (others such
as Robin can be incorporated too). In addition, we calculate $R_0$, the
radius of the smallest circle centered at the origin that encloses all
targets (see Fig.~\ref{fig:WoS1b}).

\begin{figure}[htbp]
\centering
\includegraphics[width = 0.5\textwidth]{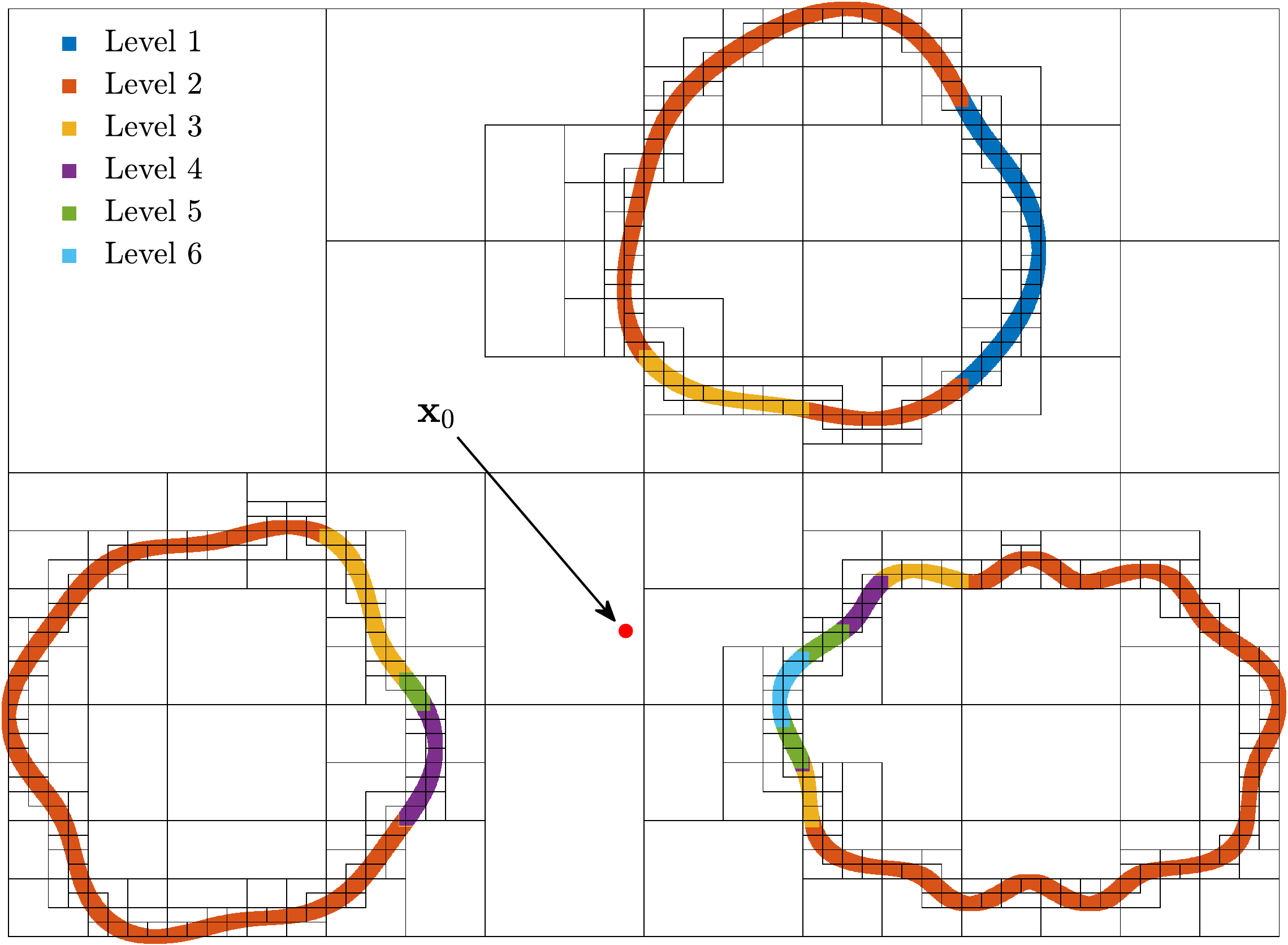}
\caption{Visualization of the tree based particle-in-cell method for
  efficient evaluation of the distance $d(\bx_0,\partial\Omega)$. The
  coloring of each point indicates the highest grid level that the
  vertex of $\partial\Omega$ is considered for the closest point
  evaluation. If the maximum refinement level is $k$, then all points in
  levels $1,2,\ldots, k-1$ are excluded from the distance calculation of
  $d(\bx_0,\partial\Omega)$.
  \label{fig:quadtree}}
\end{figure}

A frequent and potentially expensive operation is the determination of
$d(\bx_0,\partial\Omega)$, the distance of
$\bx_0\in\RR^2\setminus\Omega$ to the nearest target. A simple approach
is to calculate the distance of $\bx_0$ to each vertex of
$\partial\Omega$ and select the minimum. However, for highly refined
target geometries or numerous targets, the number of vertexes to scan
over may be prohibitive. To accommodate such scenarios, we employ a
particle-in-cell method that consists of a hierarchy of Cartesian grids
that envelop $\partial\Omega$ (Fig.~\ref{fig:quadtree}). The method
first queries the midpoints of the coarsest grid and uses simple
geometric criteria to eliminate those that cannot contain the closest
point. Queries are made of remaining subgrids at the next level of
refinement until a predefined level of refinement or a single point
remains. This results in a vastly smaller set of candidate vertexes to
calculate pointwise distances at the cost of some overheard and extends
this approach to large and complex targets sets.

With this setup in place, for each free particle $\bx_0\in
\RR^2\setminus\Omega$, we calculate the shortest distance $R =
d(\bx_0,\partial\Omega)$ and the associated projection $\textbf{p} =
\textrm{proj}_{\mathcal{P}_j} \bx_0$ where $\mathcal{P}_j$ is the line that
contains the closest edge $\partial\Omega_j$. The position of the
particle is advanced based on four projection steps
(Fig.~\ref{fig:WoS1}) described below.

\paragraph{Stage I: Radially symmetric projector}  If $R\in (R_{\textrm{min}},R_{\textrm{max}})$ such that the particle is neither too close nor too far from a target, we project to a ball of radius $R$ centered at $\bx_0$. The parameters $R_{\textrm{min}},R_{\textrm{max}}$ are associated with stages II and III respectively and defined shortly. The time duration of this projection step is determined from the solution of the radial diffusion equation with a zero Dirichlet boundary condition at $r=R$ and a Dirac initial condition specifying the particle is initially at the origin. The solution $u(r,t)$ of the parabolic equation
%\bsub\label{eq:StageI}
%\begin{equation}\label{eq:radialDiffusion2D}
\begin{align*}
\frac{\partial u}{\partial t} = D \frac{1}{r} \frac{\partial}{\partial
  r}\left(r \frac{\partial u}{\partial r} \right), \quad r\in(0,R), \quad t>0; \qquad u(R,t) = 0, \quad t>0; \qquad u(r,0) = \frac{\delta (r)}{r}, \quad r\in(0,R),
\end{align*}
%\end{equation}
gives the cumulative distribution of arrival times at $r=R$ to be
\begin{equation}\label{eq:2DArrivalTimes}
F(\tau) = 1- 2\sum_{n=0}^{\infty}\frac{ e^{-z_n^2 \tau}}{z_n J_1(z_n)}, \qquad t = \frac{D}{R^2} \tau, \qquad J_0(z_n) = 0, \quad n =0,1,2,\ldots.
\end{equation}
%\esub
This distribution is sampled by drawing a uniform number $U\in(0,1)$ and solving $F(\tau) = U$. The CDF $F(\tau)$ is efficiently sampled by precomputing values of $z_n$ and $J_1(z_n)$ and using only as many terms as is necessary to approximate $F(\tau)$ to a predetermined tolerance.

\paragraph{Stage I (alternative for convex shapes): Plane projector} This projector can be adopted when the target is strictly convex so that the entire absorbing target lies to one side of the tangent plane to the surface. After translating the projection point ${\textbf{p}}$ to the origin and rotating by the slope of the incident edge, the projection step arises from the solution to the heat equation in the upper half plane with initial location $(0,y_0)$. Combining the relevant fundamental solutions with method of images yields the density
\begin{align*}
u(x,y,t) = \frac{1}{4\pi D t} \left[e^{- \frac{x^2 + (y-y_0)^2}{4Dt}} - e^{- \frac{x^2 + (y+y_0)^2}{4Dt}} \right],
\end{align*}
with association arrival time distribution to the plane  
\begin{align*}
  \rho_T(t) = \int_{\RR} u_y(x,0,t) \,\mathrm{d}x = 
    \frac{y_0}{2\sqrt{\pi} (Dt)^{\frac32}} e^{-\frac{y_0^2}{4Dt}}.
\end{align*}
The cumulative distribution is $\int_0^t \rho_T(\tau) \,\mathrm{d}\tau =
\text{erfc}(y_0/\sqrt{4Dt})$ so that the arrival time is sampled as 
\bsub\label{eq:Plane}
\begin{align}\label{eq:PlaneTime}
t_{\ast} = \frac{1}{4D} \left[ \frac{y_0}{ \text{erfc}^{-1}(\eta)} \right]^2, \qquad \eta\in(0,1).
\end{align}
The hitting location on the tangent line is determined by the displacement $x_{\ast}$ from the projection point ${\bf p}$ which has the (Gaussian) distribution
\begin{align}\label{eq:PlaneLocation}
P_{X}(x_{\ast}) = \frac{u_y(x_{\ast},0,t_{\ast}) }{\rho_T(t_{\ast})} = \frac{1}{\sqrt{4\pi D t_{\ast}}} e^{-\frac{x_{\ast}^2}{4Dt_{\ast}} },
\end{align}
\esub
so that $x_{\ast} \sim \mathcal{N}(0,2Dt_{\ast})$.

\paragraph{Stage II: Reinsertion} It is inefficient to simulate the detailed trajectory of particles far from the absorbers. Therefore, when the distance $R = d(\bx_0,\partial\Omega)$ exceeds a threshold ($R>  R_{\textrm{max}}$), we project the particle to a smaller disc of radius $R_{\textrm{ins}} > R_0$ that encloses all the absorbers (Fig.~\ref{fig:WoS1b}). Similar reinsertion procedures have been utilized in Monte Carlo solutions of elliptic problems \cite{OKMascagni2004,DOBRAMYSL2018}. Here we must take additional care to sample both the reinsertion time and the time dependent reinsertion location correctly. The arrival distribution for a particle initially on the $x$-axis (see Appendix \ref{sec:PlanarDist}) is
\bsub\label{flux:MAINomega}
\begin{align}\label{flux:MAINomega_a}
\mathcal{J}(\tau,\theta) =   \frac{D}{R_{\textrm{ins}}^2} \left[\frac{1}{2\pi} \chi_0(\tau)  + \frac{1}{\pi}\sum_{n=1}^{\infty}\chi_n(\tau) \cos n\theta\right] , \qquad \theta \in (-\pi,\pi), \quad  t>0, \quad \tau = \frac{D}{R_{\textrm{ins}}^2} t,
\end{align}
where the coefficients are
\begin{align}
  \label{flux:MAINomega_b}
  \chi_n(\tau) = \frac{2}{\pi}\int_{0}^{\infty}\left[
    \frac{J_n(\omega)Y_n(\omega \mu)-J_n(\omega\mu)Y_n(\omega)}
    {Y_{n}^{2}(\omega)+J_{n}^{2}(\omega)}\right]\omega 
    e^{-\omega^2 \tau}\,d\omega.
\end{align}
\esub
Here $\mu := |\bx_0|/R_{\textrm{ins}}$ is the ratio of the distance
$|\bx_0|$ of the initial location from the origin to the reinsertion
radius $R_{\textrm{ins}}$. The optimal reinsertion radius is
$R_{\textrm{ins}} = R_0$ where $R_0$ is the radius of the smallest disc
enclosing all absorbers as shown in Fig.~\ref{fig:WoS1b}. However, many
computational efficiencies are gained by sampling \eqref{flux:MAINomega}
for a fixed value of $\mu$ - in practice we take $\mu = 60$,
$R_{\textrm{max}} = \mu R_0$ and reinsert to the disc of radius
$R_{\textrm{ins}} = |\bx_0|/\mu$. By fixing $\mu$, the integrands of
\eqref{flux:MAINomega_b} can be tabulated over a range of $\omega$
values for efficient quadrature. 

\begin{figure}[htbp]
\centering
\subfigure[Arrival time distribution.]{\includegraphics[width = 0.4\textwidth]{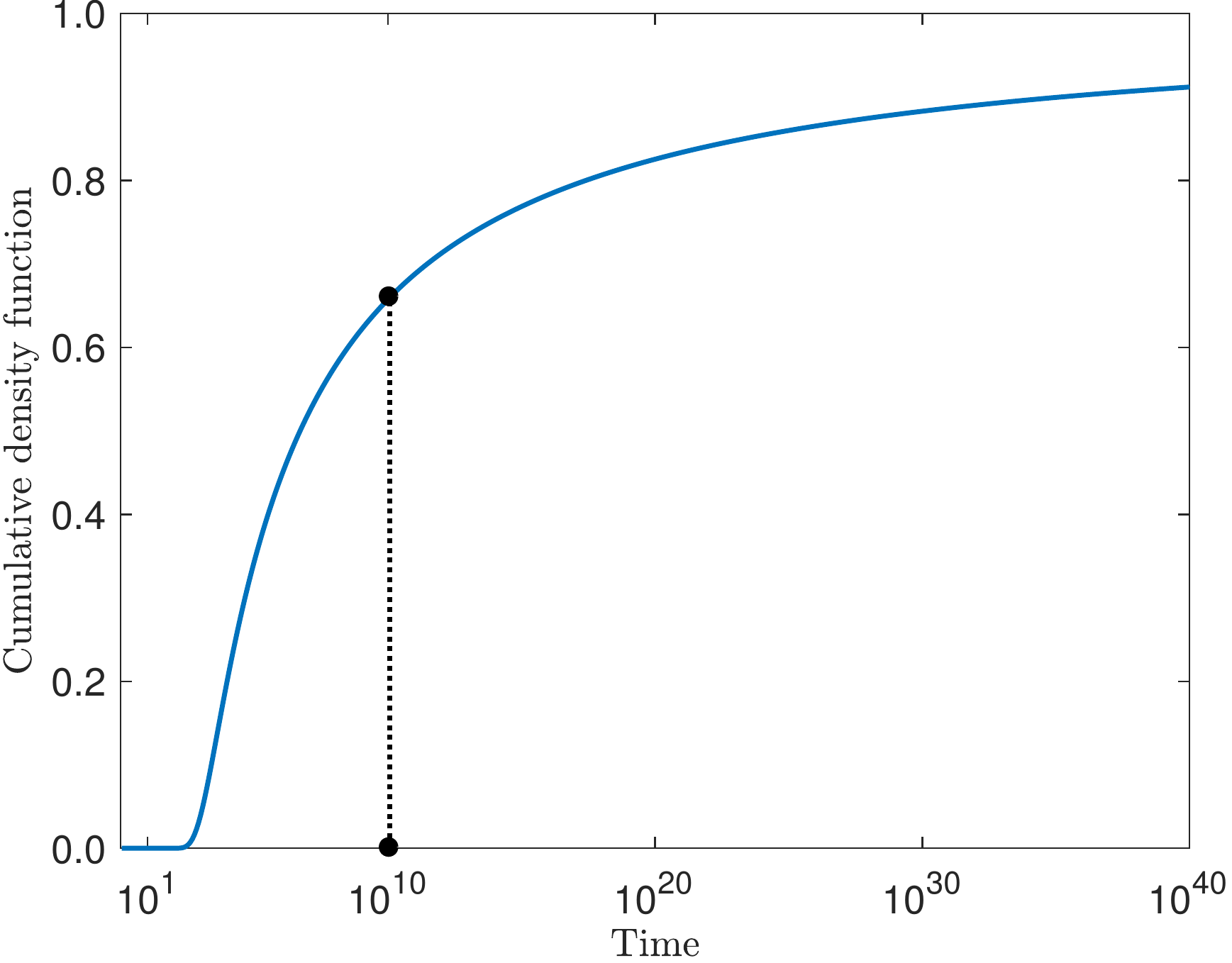}\label{fig:2DarrivalCDF_a}}\hspace{0.5in}
\subfigure[Arrival angle distribution for several $t_{\ast}$ values.]{\includegraphics[width = 0.39\textwidth]{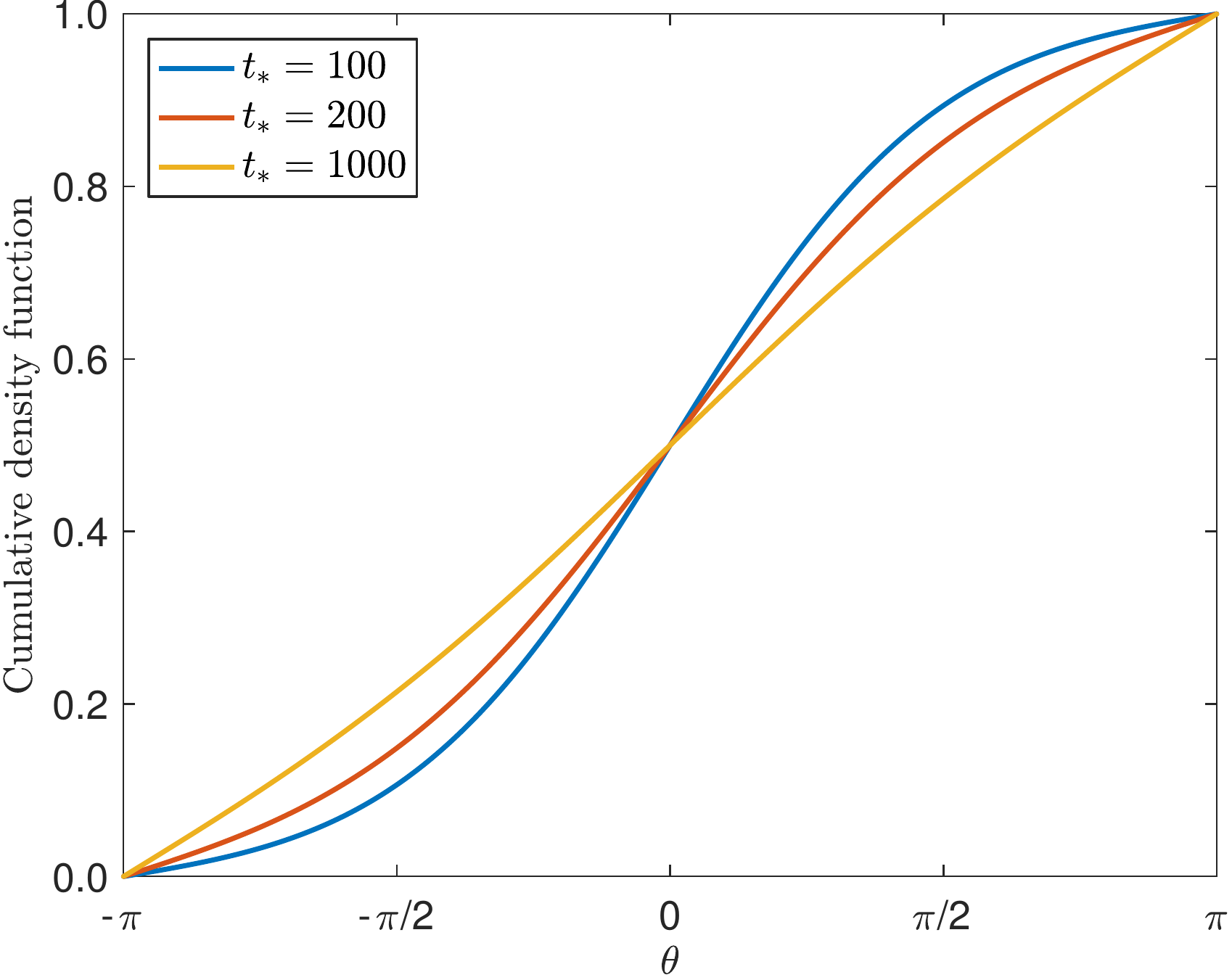}\label{fig:2DarrivalCDF_b}}
\caption{CDFs of arrival time and location for the reinsertion step
  when $D=1$, $\mu=60$. Left: The reinsertion time distribution
  $F_{\theta}(\tau) = \int_0^{\tau} \chi_0(\eta)\,d\eta$ given by
  \eqref{eq:dist2Dreinsertion} (dashed line indicates the cutoff
  $\tau =10^{10}$ after which the fitted density
  \eqref{eq:fittedDensity} is used). Right: The arrival location CDF
  \eqref{eq:samp:Th} for several arrival times $\tau_{\ast}$. The distribution is largely uniform for large $\tau_{\ast}$. \label{fig:2DarrivalCDF}}
\end{figure}

The first step in the sampling of \eqref{flux:MAINomega} is to determine
the arrival time density $\int_{\theta=0}^{2\pi} \mathcal{J}(
\tau,\theta)\,d\theta = \frac{D}{R_{\textrm{ins}}^2}\chi_0( \tau) $ with
associated CDF $F_T(t) = \int_0^{ \tau} \chi_0(\eta)\,d\eta$ where $\tau
= \frac{D}{R_{\textrm{ins}}^2}t$. The arrival time is sampled first
since the location will be dependent on this value---for shorter times
the arrival location is more tightly focussed around the initial
location while for larger arrival times, the insertion location has a
weaker dependence on the start location and approaches a uniform
distribution (see Fig.~\ref{fig:2DarrivalCDF_b}). For smaller values of
$\tau$\, (in practice $\tau<10^{10}$), the values of the integrand \eqref{flux:MAINomega_b} and
the associated CDF are tabulated over a range of $\omega$ values for
rapid quadrature. The sampling of $\chi_0(\tau)$ can be quite delicate
for large $\tau$ due to the slow convergence of the integral. To see
this, consider that for $\tau\gg1$ the main contribution to the integral
$\chi_0(\tau)$ is when $\omega^2 \tau = \bigoh(1)$ or $\omega\ll1$. In
this regime we have that,
\begin{align}\label{eq:ApproxLarget}
\chi_0(\tau) \sim \frac{4 \log \mu}{\pi^2} \int_{0}^{\infty} \frac{w}{1+
  \frac{4}{\pi^2}(\gamma_e + \log(w/2))^2} e^{-\omega^2 \tau} \, d\omega
  = \bigoh\left(\frac{1}{\tau|\log \tau|^2}\right), \qquad \tau \to \infty.
\end{align}
When $\tau\gg1$ (in practice $\tau >10^{10}$), we use the limiting form \eqref{eq:ApproxLarget} to posit the following explicit form of the density
\begin{align*}
  \chi_0(\tau) = \frac{4\log \mu}{\pi^2}
  \frac{1}{\tau(a_1 + a_2 \log \tau + a_3 \log \tau ^2)} + 
  \bigoh\left(\frac{1}{\tau^2}\right), \qquad \tau\gg1,
\end{align*}
for constants $a_1, a_2, a_3$ determined from fitting. This gives the exact cumulative density function for $\tau\gg1$
\begin{align}\label{eq:fittedDensity}
F_T(t) = 1- \int_{\tau}^{\infty} \chi_0(\eta)\,d\eta = 1 - \frac{4\log
  \mu}{\pi^2\sqrt{4a_1 a_3 - a_2^2}} \left( \pi - 2\tan^{-1}
  \left[\frac{a_2 + 2a_3 \log \tau}{\sqrt{4a_1 a_3 - a_2^2}}\right]
  \right), \quad \tau = \frac{D}{R_{\textrm{ins}}^2} t,
\end{align}
where for $\mu = 60$, we obtain from fitting the constants
\begin{align*}
a_1 = 1.4670, \qquad a_2 = 0.3102, \qquad a_3 = 0.2029.
\end{align*}
For a particular arrival time realization $\tau = \tau_{\ast}$, the angular location $\theta$ of reinsertion satisfies
\begin{equation}\label{eq:samp:Th}
F_{\theta}(\theta; \tau_{\ast}) = \int_{-\pi}^{\theta}
  \frac{\mathcal{J}(\tau_{\ast},\eta)}{\chi_0(\tau_{\ast})}\,d\eta = \frac{\theta+\pi}{2\pi} + \frac{1}{\pi}\sum_{n=1}^{\infty} \frac{\sin (n \theta)}{n} \frac{\chi_n( \tau_{\ast})}{\chi_0(\tau_{\ast})}, \qquad \theta\in(-\pi,\pi).
\end{equation}
To sample from \eqref{eq:samp:Th}, a uniform number $U\in(0,1)$ is drawn and the equation $F_{\theta}(\theta;\tau_{\ast}) = U$ solved for $\theta$. In practice we use an adaptive procedure where we retain only the terms satisfying $|\chi_n(\tau_{\ast})/(n\chi_0(\tau_{\ast}))| > \tau_{\textrm{tol}}$ in the summation of \eqref{eq:samp:Th}. For a large proportion of realizations, the arrival time $\tau_{\ast}$ is sufficiently large (in practice $\tau_{\ast}>10^{10}$) so that the first term is negligible and the arrival location is uniformly distributed on the disc. In Fig.~\ref{fig:2DarrivalCDF}, we plot the CDFs of arrival time distribution and arrival location distribution for $\mu=60$. This method permits rapid and accurate sampling of the reinsertion step. 

\paragraph{Stage III: Square projector} 
If $R<R_{\textrm{min}}$, then the particle is close enough to determine if contact occurs. By \lq\lq close enough\rq\rq, we mean that the projection $\bf{p}$ lies within the edge segment so that a square of side length $2R$ centered at $\bx_0$ lies entirely within the target edge  (cf.~Fig.~\ref{fig:WoS1c}). This gives explicitly that $R_{\textrm{min}} = \min(d_1,d_2)$ where $d_1$, $d_2$ are the distances between $\bf{p}$ and the edge vertexes (cf.~Fig.~\ref{fig:WoS1c}). The projection step is then determined from the solution to the parabolic equation on the square $\mathcal{S} = [-R, R]^2$
\bsub\label{eq:SquareDiffusion2D}
\begin{gather}
 \label{eq:SquareDiffusion2D_a} \frac{\partial u}{\partial t} = D \left(
  \frac{\partial^2 u}{\partial x^2} + \frac{\partial^2 u}{\partial
  y^2}\right), \quad \mathbf{x} \in \mathcal{S}, \quad t > 0; \quad u(\mathbf{x}, t) = 0, \quad \mathbf{x} \in \partial \mathcal{S}, \quad t>0; \\[5pt] \quad u(\mathbf{x}, 0) = \delta(x)\delta(y), \quad \mathbf{x} \in \mathcal{S}, \quad t= 0.
\end{gather}
\esub
The separable solution to \eqref{eq:SquareDiffusion2D} yields the CDF of first arrival times to the square edge $\partial\mathcal{S}$
%\bsub\label{eq:2DSquareArrival}
%\begin{equation}\label{eq:2DSquareArrivalTimes}
\begin{align*}
P_T(\tau) = \int_0^{\tau} \rho_T(\eta)\,d\eta = 1-\frac{32}{\pi^2}  \sum_{l=0}^{\infty} \sum_{k=0}^{\infty}  \frac{2k+1}{2l+1} \frac{(-1)^{l+k}}{(2k+1)^2 + (2l+1)^2} e^{-((2l+1)^2 + (2k+1)^2) \pi^2 \tau}. \quad \tau = \frac{D}{R^2}t. 
\end{align*}
% \end{equation}
This distribution is sampled by drawing a uniform number $U\in(0,1)$ and solving $P_T(\tau^{\ast}) = U$. Each side of the square has an equal probability $1/4$ of being hit and the arrival location is sampled from the density
%\begin{equation}\label{eq:2DSquareArrivalLocation}
  \begin{align*}
\rho_X(x) = \frac{4u_y(x,0,\tau^{\ast})}{\rho_T(\tau^{\ast})} = \frac{\pi}{2}\, \ds\frac{\ds\sum_{l=0}^{\infty} \ds\sum_{k=0}^{\infty}  (2l+1) (-1)^{l+k}  \, e^{-((2l+1)^2 + (2k+1)^2) \pi^2 \tau^{\ast}}\sin[(2k+1)\pi x]}{\ds\sum_{l=0}^{\infty} \ds\sum_{k=0}^{\infty}  \frac{2k+1}{2l+1} (-1)^{l+k} e^{-((2l+1)^2 + (2k+1)^2) \pi^2 \tau^{\ast}}}.
\end{align*}
%\end{equation}
%\esub
In practice, we precalculate a large number of pairs $\{(\tau_j,x_j)\}$ to draw from when needed during simulation.
\paragraph{Stage IV: Reflection step} 

In the scenario that the particle hits a reflecting portion of the target surface, it is then projected back into the bulk onto a semi-circle of radius $R = \min(d_1,d_2)$, corresponding to the distance to the nearest vertex. In practice, we avoid rounding errors by setting $R = \max(\epsilon,\min(d_1,d_2))$ where $\epsilon$ is a small number comparable to machine precision. In the reflecting boundary condition scenario, the projection step is identical to that of stage I with uniform location and arrival time sampled from \eqref{eq:2DArrivalTimes}. For highly convoluted geometries, it may be that this semi-circle intersects with distal elements of the target. This scenario can be accounted for by calculating the shortest distance $d$ to other segments of the boundary and setting $R = \max(\epsilon,\min(d_1,d_2,d))$.

\subsection{Rudimentary convergence analysis in half plane case}

\begin{figure}[H]
\centering
\subfigure[KMC step]{\includegraphics[clip, width = 0.425\textwidth]{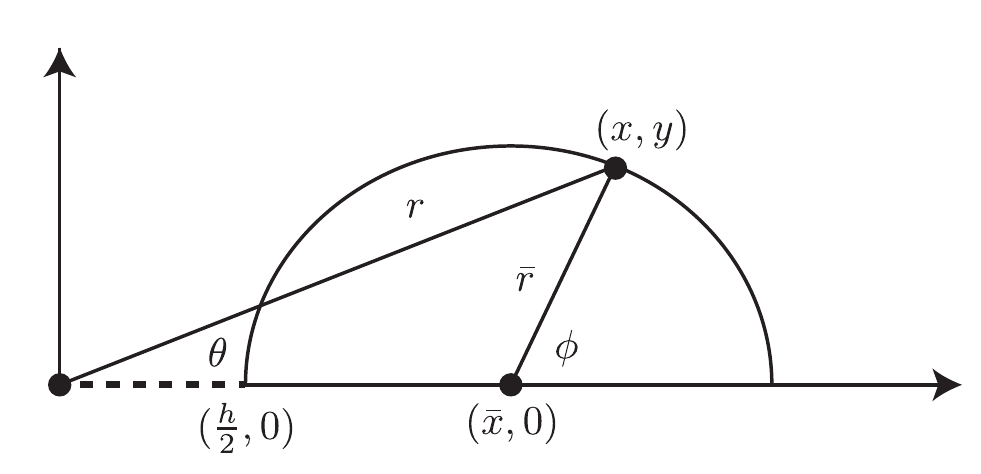}}
\qquad
\subfigure[Reinsertion]{\includegraphics[clip, width = 0.45\textwidth]{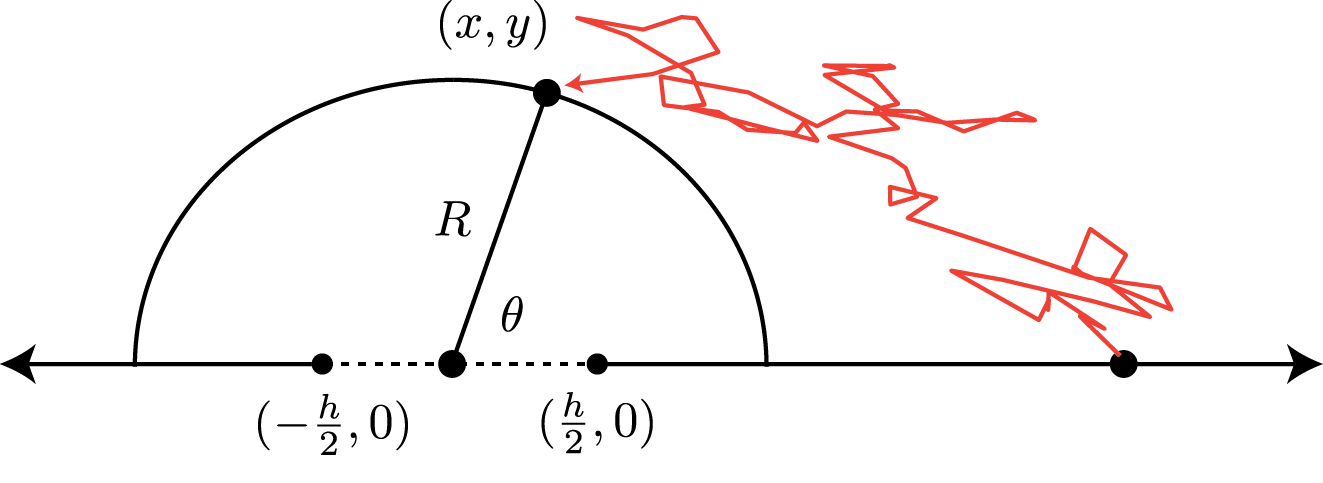}}
\caption{A schematic showing the simplified KMC solution in the upper half place. (a) the projection step from the real line back into the bulk. (b) The reinsertion step to keep the particle in the vicinity of the absorbing set. \label{Fig:reinsertion_schematic}}
\end{figure}

In this section we present some analysis of the convergence properties of the KMC approach in the simplified scenario of diffusion in the upper half plane with capture in the window $\{(x,0) \, \mid \,  |x|<h/2\}$ with particular emphasis on the role of reinsertion. Specifically, we solve the equation
\bsub\label{eq:half_plane_full}
\begin{gather}
\label{eq:half_plane_full_a}    \frac{\partial u}{\partial t} =
  \frac{\partial^2 u}{\partial x^2} + \frac{\partial^2 u}{\partial y^2},
  \quad x\in\RR, \quad y >0, \quad t>0; \\[5pt]
\label{eq:half_plane_full_b}     u = 1, \quad y = 0, \quad|x| < \frac{h}{2}, \qquad u_y = 0,\quad y = 0, \quad |x| \geq \frac{h}{2};\\[5pt]
\label{eq:half_plane_full_c}  u = \delta(x-x_0)\delta(y-y_0), \quad
  x\in\RR, \quad y >0, \quad t=0.
\end{gather}
\esub
with a simplified KMC method composed of Stages 1(alternative) and 3.
From an initial point $\bx_0=(x_0,y_0)$, this algorithm results in a
sequence of points $(\bx_1,\bx_2,\ldots, \bx_n)\in \RR^{2+}$ which will eventually alight on the absorbing portion.  For an ensemble of particles, we denote $a_n$ to be the fraction free after $n$ iterations so
\[
a_{n+1} =  (1-{p_n}) a_{n}, \qquad a_0 = 1,
\]
where ${p_n}$ is the probability of capture at the $n^{th}$ iteration. To investigate $p_n$, we first consider the splitting problem $p(x,y)$ for the probability that a particle starting at $(x,y)$ first contacts the plane $y=0$ on the absorbing window. This satisfies
\bsub\label{split}
\begin{gather}
\label{split_a}    \frac{\partial^2 p}{\partial x^2} +  \frac{\partial^2
  p}{\partial y^2} = 0, \quad x\in\RR, \quad y>0;\\[5pt]
\label{split_b}       p = 1, \qquad y = 0, \quad |x|< h/2, \qquad p =0, \qquad y =0, \quad |x|\geq h/2,
\end{gather}
and admits the solution
\begin{equation}\label{split_c}   
    p(x,y) = \frac{1}{\pi}\left[ \tan^{-1}\left(\frac{x+\frac{h}{2}}{y} \right) - \tan^{-1}\left(\frac{x-\frac{h}{2}}{y} \right)  \right].
\end{equation}
\esub
We consider that the current location $\bx_n = (x_n,y_n)$ has arisen from a projection step (stage III). Without loss of generality, we assume the previous contact $(\bar{x}_n,0)$ with the $y_n=0$ plane is such that $\bar{x}_n> \frac{h}{2}$. We may then parameterize (cf.~Fig.~\ref{Fig:reinsertion_schematic}) the point $\bx_n = (x_n,y_n)$ as 
\begin{align*}
x_n = \frac{h}{2} + \bar{r}_n(1+\cos\phi), \qquad y_n = \bar{r}_n\sin\phi, \qquad \phi \in(0,\pi), \qquad \bar{x}\in(h/2,\infty).
\end{align*}

It follows from the splitting probability \eqref{split_c} that
\begin{align}
\nonumber p(x_n,y_n) &= \frac{1}{\pi}\left[ \tan^{-1} \left( \frac{x_n-\frac{h}{2}}{y_n} + \frac{h}{y_n} \right) - \tan^{-1} \left( \frac{x_n-\frac{h}{2}}{y_n}\right)\right]\\[5pt]
\nonumber {} & = \frac{1}{\pi} \left[  \tan^{-1} \left( \cot \frac{\phi}{2} + \frac{h}{\bar{r}_n} \csc \phi \right) - \tan^{-1} \left( \cot \frac{\phi}{2}\right) \right]\\[5pt]
\nonumber {} & = \frac{1}{\pi} \left[  \tan^{-1} \left( \tan \left(\frac{\pi}{2} - \frac{\phi}{2}\right) + \frac{h}{\bar{r}_n} \csc \phi \right) - \left(\frac{\pi}{2} - \frac{\phi}{2}\right) \right],
\end{align}
where the angles $\phi\in(0,\pi)$ are distributed uniformly. Applying the change of variables $\eta = \frac{\pi-\phi}{2}$ yields
\begin{equation*}
    p = \frac{1}{\pi} \left[  \tan^{-1} \left( \tan \eta + \alpha_n \csc 2\eta \right) - \eta \right], \qquad \alpha_n = \frac{h}{\bar{r}_n}.
\end{equation*}
We now define the average probability $p_n =
\frac{1}{\pi}\int_{\phi=0}^{\pi} p \, d\phi =
\frac{2}{\pi}\int_{\eta=0}^{\pi/2} p \, d\eta$ and find that
\begin{align}\label{eqn:pn}
p_n = \frac{2}{\pi^2}\int_{\eta = 0}^{\frac{\pi}{2}} \left[ \tan^{-1}
  \left( \tan \eta + \alpha_n \csc 2\eta \right) - \eta \right] \, d\eta, \qquad \alpha_n = \frac{h}{\bar{r}_n}.
\end{align}

Without reinsertion, many trajectories will yield very small values for the parameter $\alpha_n$ as $\bar{r}_n$ can attain very large values. To gain further insight into this scenario, we consider the limiting case of \eqref{eqn:pn} as $\alpha\to0$. 

\subsection{Asymptotic analysis of splitting probability}

Here we develop an asymptotic approximation for \eqref{eqn:pn} in the limit as $\alpha\to0$ (subscript $n$ dropped for convenience). The integral features global contributions and local contributions near $\eta=0$. To delineate between these contributions, we define the small parameter $\delta>0$ such that $\alpha \ll \delta \ll 1$. Then we have that
\begin{align}\label{eqn:pn_alpha}
    p_n(\alpha) = \underbrace{\frac{2}{\pi^2}\int_{\eta = 0}^{\delta}
    \left[ \tan^{-1} \left( \tan \eta + \alpha \csc 2\eta \right) - \eta
    \right] \, d\eta}_{A_1} +
    \underbrace{\frac{2}{\pi^2}\int_{\eta = \delta}^{\frac{\pi}{2}}
    \left[ \tan^{-1} \left( \tan \eta + \alpha \csc 2\eta \right) - \eta
    \right] \, d\eta}_{A_2}.
\end{align}
where $A_1$ and $A_2$ will be considered separately and combined so that their sum is independent of $\delta$.

\paragraph{Evaluation of $A_1$} In this region we apply the
approximation $\tan^{-1}(x + y) = \tan^{-1}(x) + y/(1+x^2) + \bigoh(y^2)$ for $y\to0$ with $x = \alpha \csc 2\eta$ and $y = \tan \eta$. Then we have that
\begin{equation}
 \nonumber A_1 := \frac{2}{\pi^2}\int_{0}^{\delta} \left[ \tan^{-1}
  \left( \tan \eta + \alpha \csc 2\eta \right) - \eta \right]\, d\eta \sim
  \frac{2}{\pi^2} \int_{0}^{\delta} \left[ \tan^{-1} (\alpha\csc 2\eta)
  + \frac{\tan \eta}{1+ \alpha^2 \csc^2(2\eta)} - \eta\right]\, d\eta.
\end{equation}
Applying small argument approximations for $\eta\ll1$, we have that
\begin{align}
\nonumber A_1 &\sim \frac{2}{\pi^2}\int_{0}^{\delta}\left(
  \tan^{-1}\left(\frac{\alpha}{2\eta}\right)  - \frac{\alpha^2
  \eta}{4\eta^2 + \alpha^2}\right) \, d \eta = \frac{\alpha}{\pi^2}
  \int_0^{\tfrac{2\delta}{\alpha}}\left( \tan^{-1} \frac{1}{z} -
  \frac{\alpha}{2} \frac{z}{z^2+1}\right) \, dz, \qquad (\eta = \alpha z/2)\\[5pt]
\nonumber {} & =  \frac{\alpha}{\pi^2} \left[\frac{1}{2} \log(1+ z^2) +
  z \tan^{-1}\frac{1}{z} - \frac{\alpha}{4}\log(1+ z^2)
  \right]_0^{\frac{2\delta}{\alpha}} = \frac{\alpha}{\pi^2} \left[
    \left(\frac{1}{2} - \frac{\alpha}{4}\right) \log\left(1 +
    \frac{4\delta^2}{\alpha^2}\right) + \frac{2\delta}{\alpha} \tan^{-1}
    \frac{\alpha}{2\delta}  \right].
\end{align}
The parameter $\delta$ is defined so that $\alpha \ll \delta$ and therefore $\frac{\delta}{\alpha}\gg1$. Applying Taylor we find that
\begin{align}\label{eqn:A1}
A_1 = \frac{\alpha}{\pi^2} \left( \log \frac{2\delta}{\alpha} + 1
  \right) + \bigoh\left(\frac{\alpha^3}{\delta^2},
  \alpha^2\log\left(\frac{\delta}{\eps}\right)\right).
\end{align}

\paragraph{Evaluation of $A_2$} In the integral $A_2$ we
have that $\alpha\ll\delta$ and apply the approximation $\tan^{-1}(x +
y) = \tan^{-1}(x) + y/(1+x^2) + \bigoh(y^2)$ for $y\to0$ with $x = \tan \eta$ and $y = \alpha \csc 2\eta$. It then follows that
\begin{align}
\nonumber A_2 &:= \frac{2}{\pi^2}\int_{\eta = \delta}^{\frac{\pi}{2}}
  \left[ \tan^{-1} \left( \tan \eta + \alpha \csc 2\eta \right) - \eta
  \right] \, d\eta \sim \frac{2}{\pi^2}\int_{\delta}^{\frac{\pi}{2}}
  \left[\tan^{-1} (\tan\eta) + \alpha \frac{\csc 2\eta}{1 + \tan^2\eta}
  - \eta \right] \, d\eta\\[5pt]
\nonumber & = \frac{2\alpha}{\pi^2}\int_{\delta}^{\frac{\pi}{2}}
  \frac{\csc 2\eta}{1 + \tan^2\eta} \, d\eta = \frac{\alpha}{\pi^2}
  \int_{\delta}^{\frac{\pi}{2}} \cot\eta \, d\eta = \frac{\alpha}{\pi^2}
  \left[\log\sin\eta \right]_{\delta}^{\frac{\pi}{2}}\\[5pt]
& = -\frac{\alpha}{\pi^2} \log\delta. \label{eqn:A2}
\end{align}

To finalize the approximation of \eqref{eqn:pn_alpha}, we combine expressions \eqref{eqn:A1} and \eqref{eqn:A2} and reintroduce $\alpha_n$ reflecting that this parameter changes over each iteration. This yields that 
\begin{equation}\label{eqn:pn_asy}
    p_n \sim A_1 + A_2 = \frac{\alpha_n}{\pi^2} \left(\log \frac{2}{\alpha_n} + 1 \right), \qquad \alpha_n \to 0, \qquad \alpha_n = \frac{h}{\bar{r}_n}.
\end{equation}
Hence we see that while the probability of capture is positive at each iteration, it can become arbitrarily small as $|\bar{x}_n| =  \frac{h}{2} + \bar{r}_n \to\infty$. This results in an algorithm with polynomial convergence rate (see Fig.~\ref{fig:reinsertion_convergence_a}).

\begin{figure}
    \centering
\subfigure[Capture probability without reinsertion.]{\includegraphics[width= 0.45\textwidth]{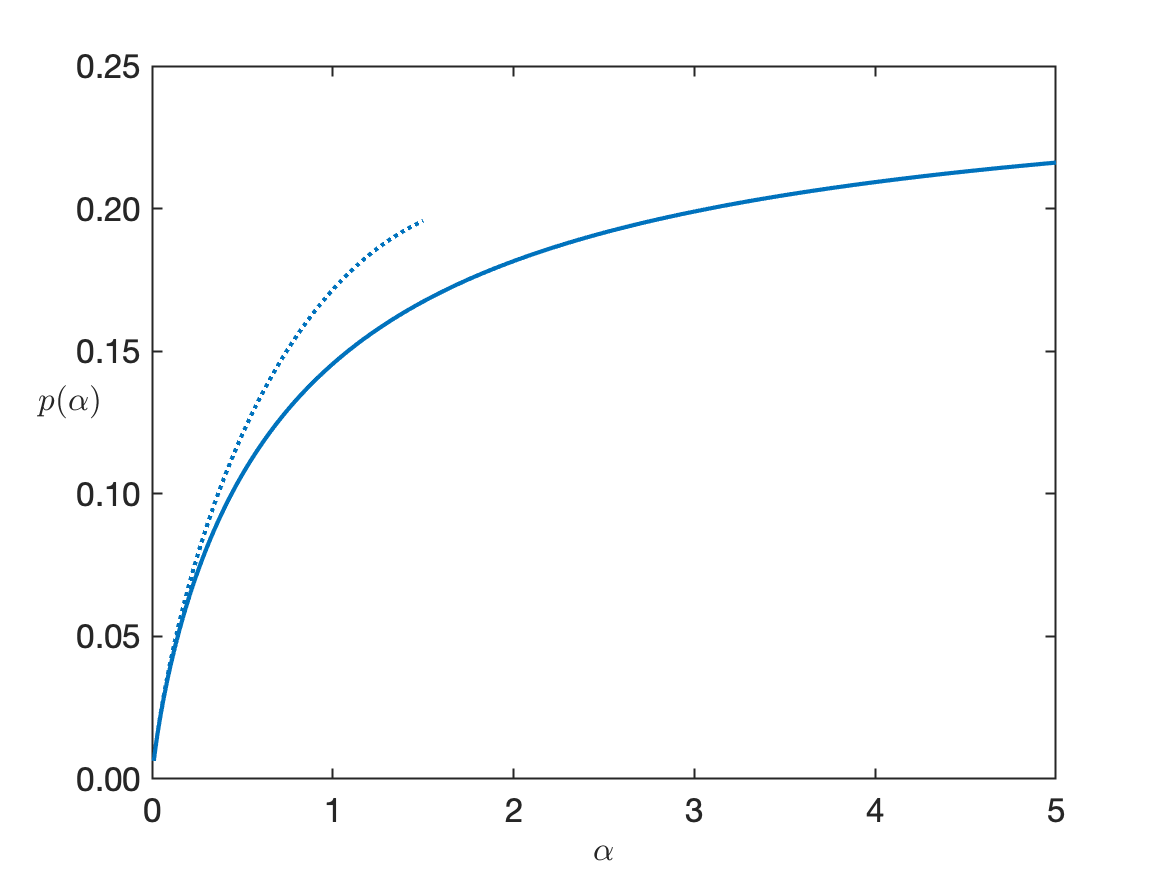}}
\subfigure[Capture probability from reinsertion disk.]{\includegraphics[width= 0.45\textwidth]{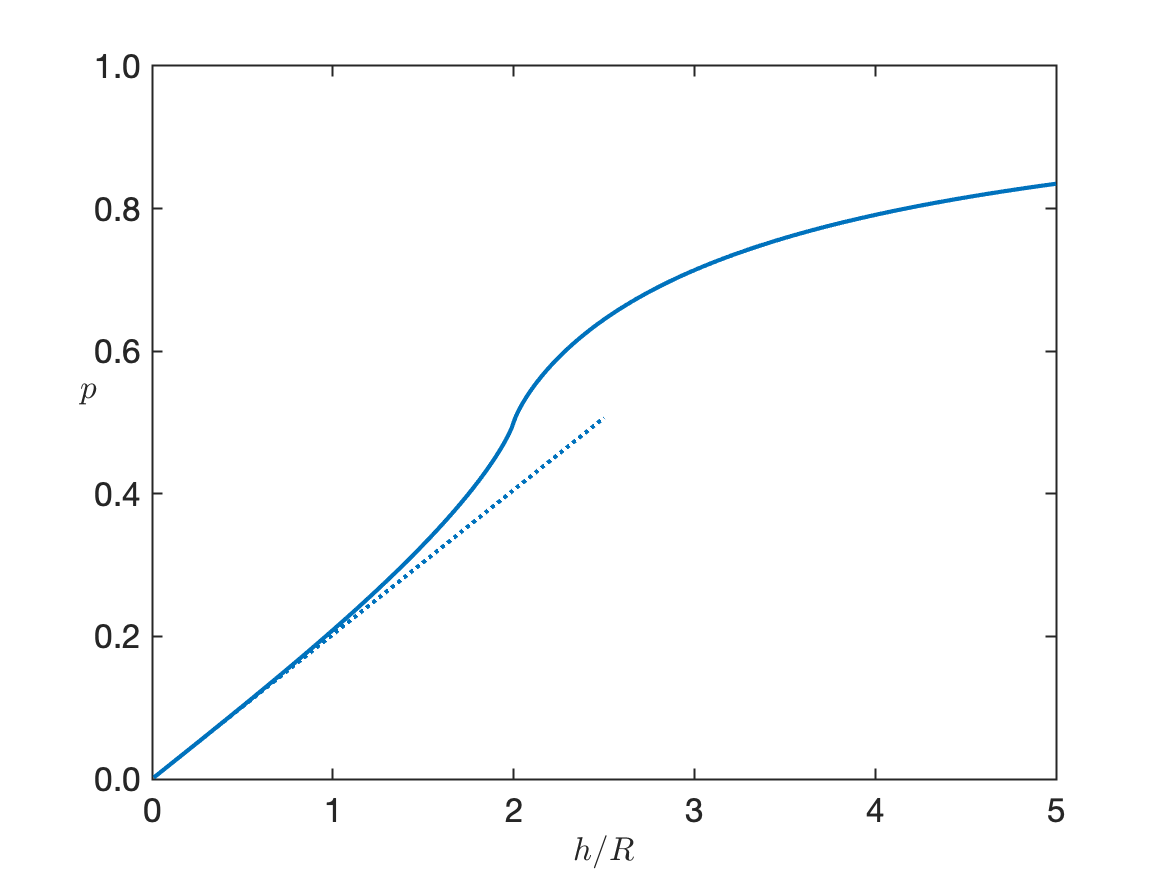}}
    \caption{(Left figure) Comparison of the capture probability $p_n(\alpha)$ (solid curve) \eqref{eqn:pn} and the small arguement approximation \eqref{eqn:pn_asy} (dashed curve) without reinsertion. We note that $p\to 1/4$ for $\alpha_n = h/\bar{r}_n \to \infty$. (Right figure) Here we show the probability of capture with reinsertion to a disk of radius $R$. }
    \label{fig:my_label}
\end{figure}

\subsection{Reinsertion analysis}

The aim of reinsertion is to reestablish exponential convergence rate in the algorithm by limiting the maximum value of $\bar{x}_n$ and hence promoting faster capture.

Reinsertion projects wayward particles back to a smaller disc of radius $R$ that encloses all targets and therefore omits simulating  trajectories far from the capture regions. When reinserting from sufficiently distant points, the placement on the disc is largely uniform with $(x,y) = R(\cos\theta,\sin\theta)$ with small corrections given by \eqref{eq:samp:Th}. Points on this disk have average probability of capture (see \eqref{split})
\begin{align}
\nonumber p_n &=
   \frac{1}{\pi}\int_{0}^{\pi}p(R\cos\theta,R\sin\theta)\,d\theta =
  \frac{1}{\pi^2} \int_{0}^{\pi} \left[  \tan^{-1} \big( \cot \theta +
  \frac{h}{2R} \csc \theta \big) - \tan^{-1} \big( \cot \theta -
  \frac{h}{2R} \csc \theta \big) \right]\, d\theta\\[5pt]
\label{prob_Rins} &\sim \frac{2}{\pi^2} \frac{h}{R}, \quad \mbox{as} \quad \frac{h}{R} \to 0.  
\end{align}
In an ensemble of particles, some will be reinserted to the disk $x^2+y^2 = R^2$ while others remain inside it. Those inside have greater probability of capture at the next step, therefore the quantity \eqref{prob_Rins} reflects a lower bound on the likelihood of capture. Hence we see that the probability of capture after $n$ stages has bound $p_n \gtrapprox \frac{2}{\pi^2} \frac{h}{R}$. The key observation here is that the probability of capture at each iteration is now bounded below by a constant defined in term of the geometric parameter $h$ and the reinsertion radius $R$.  This ensures an exponential convergence rate of the algorithm with slower rates associated with smaller targets ($h\ll1$) and a faster rate associated with a smaller reinsertion radius $R$.

\begin{figure}[htbp]
\centering
\subfigure[Convergence without reinsertion]{\includegraphics[clip, width
  = 0.45\textwidth]{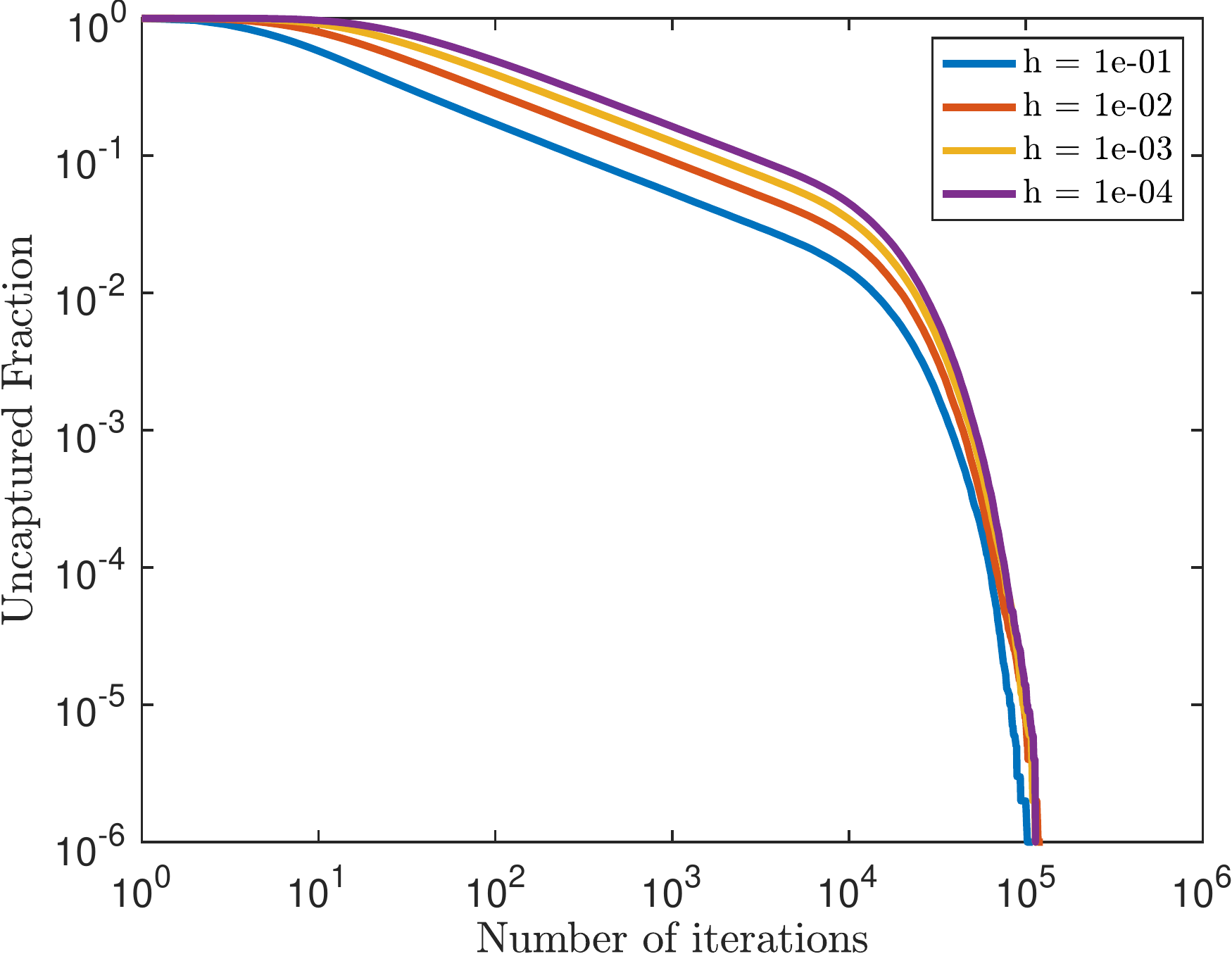}\label{fig:reinsertion_convergence_a}}
\qquad
\subfigure[Convergence with reinsertion]{\includegraphics[clip, width =
  0.45\textwidth]{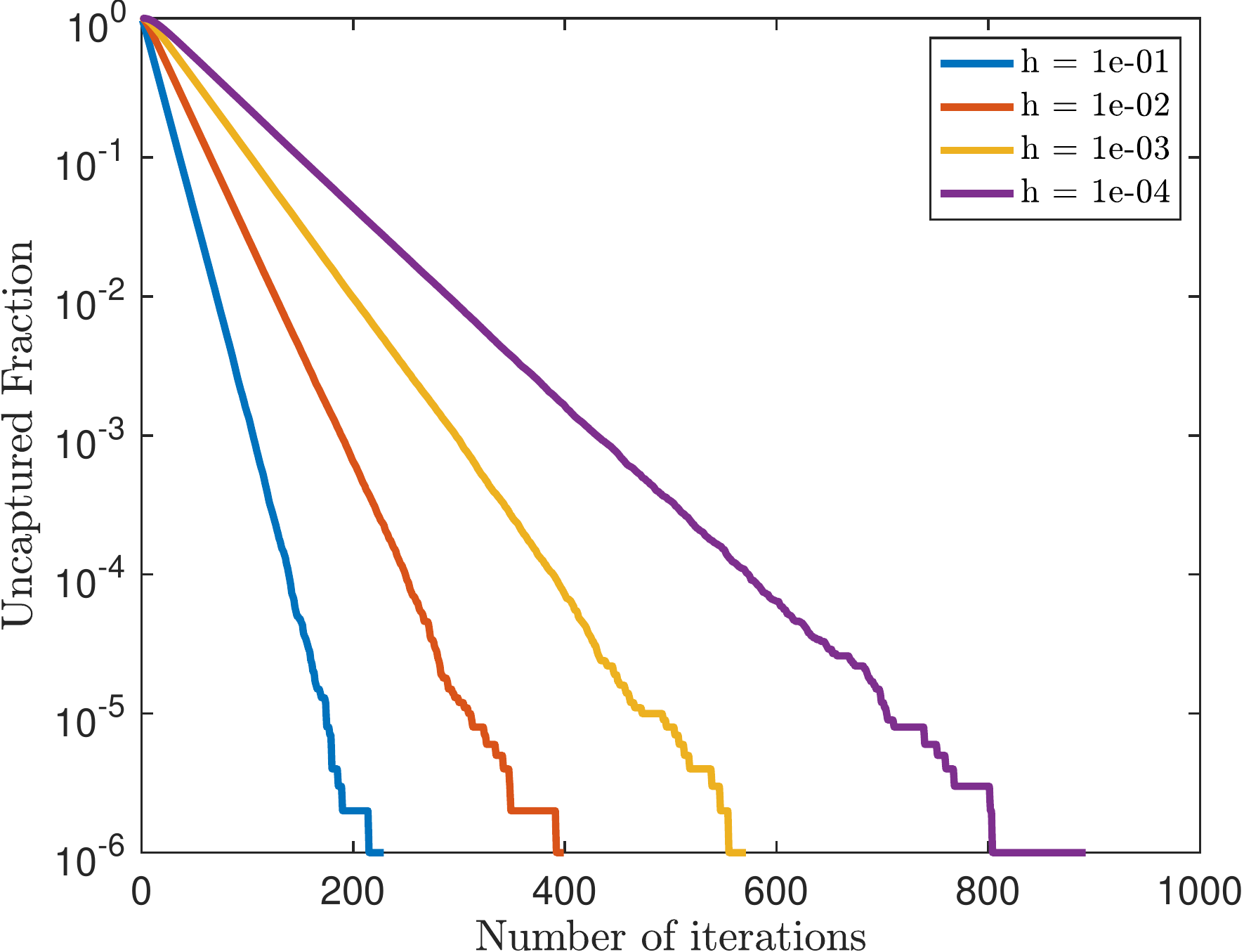}\label{fig:reinsertion_convergence_b}}
\caption{Convergence statistics of the KMC method for the simplified problem \eqref{eq:half_plane_full}. Left: In the absence of reinsertion, a polynomial rate of convergence is attained shown by the linear behavior on the log-log plot. Right: Reinsertion recovers exponential convergence as seen in the linear behavior on the log plot. Simulations based on ensemble of $10^6$ particles. \label{fig:reinsertion_convergence}}
\end{figure}

As an exposition of this analysis, we show in Fig.~\ref{fig:reinsertion_convergence} the convergence of the KMC method on the simplified half plane problem \eqref{eq:half_plane_full} with and without reinsertion. In the absence of reinsertion, polynomial convergence is attained as shown by linear behavior on a log-log plot. When reinsertion is implemented (to radius $R = 1$), we observe exponential convergence shown by linear behavior on a log plot. This demonstrates the key role of reinsertion in attaining exponential convergence of the KMC method.

\section{Results}\label{sec:results}

In this section we show a variety of examples to demonstrate the ability of these methodologies for approximating full arrival time distributions to complex absorbing sets. In our simulations, we have used a common diffusivity of $D=1$. Arrival times from the KMC method are translated with $t\to t+1$ so that $t=0$ is mapped to $10^0$ in log space. \texttt{MATLAB}'s histogram function is then applied with the \lq\lq probability\rq\rq\ normalization option. The hybrid approaches refer to solving the Laplace transform using either the expansion asymptotic \eqref{sysmain} or the boundary integral method (BIM) described in Sec.~\ref{sec:BIE}. This is followed by numerical inversion of the transform equation as described in Sec.~\ref{sec:invLT}.

\subsection{Planar results}\label{sec:resultsPlanar}

Here we use three examples to validate the numerical KMC method and corroborate with both the hybrid approaches (asymptotic and boundary integral). The first example is a simple one target scenario in which the closed form solution \eqref{eqn:timepdf} is available for comparison. The remaining two examples show the efficacy of the method on more complex absorbing sets consisting of multiple targets of varying radii. The parameter values for the three examples are
\begin{align*}
(\mbox{One target}) & \quad \bx_1 = (0,0), \quad r_1 = 0.05, \quad \bx_0 = (5,0).\\[5pt]
(\mbox{Three targets})& \quad \bx_1 = (3,3), \quad \bx_2 = (8,8), \quad \bx_3 = (10,10), \qquad r_{1}= \frac18, \quad r_{2} = \frac16, \quad r_{3} = \frac13, \quad \bx_0 = (0,0).\\[5pt]
(\mbox{Six targets})& \quad \bx_1 = (-3,0), \quad  \bx_2 =(0,-2), \quad  \bx_3 =(\sqrt{2},-\sqrt{2}),\quad \bx_4 =(2,0),\\[5pt]
\phantom{(\mbox{Six targets})}& \quad \bx_5 =(\sqrt{2},\sqrt{2}),\quad \bx_6 =(0,2), \qquad
r_1  = 0.275, \quad r_{2} = \ldots =r_{6} = 0.02, \quad \bx_0 = (0,0).
\end{align*}

\begin{figure}[htbp]
 \vspace{-1\baselineskip}% Squeeze space before.
  \centering
   \subfigure[One target example]{\includegraphics[width=0.325\linewidth]{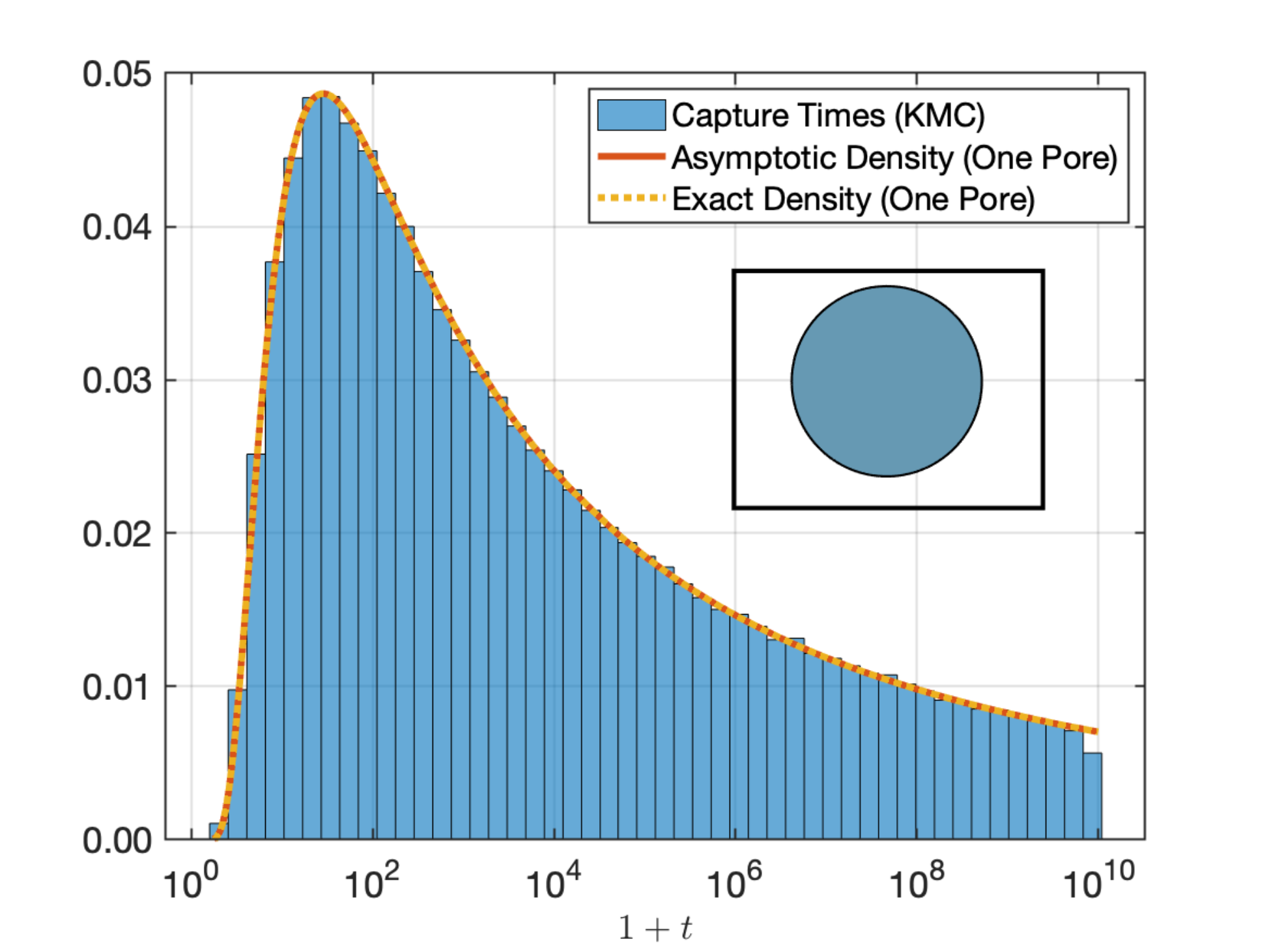}}
   \subfigure[Three target example]{\includegraphics[width=0.325\linewidth]{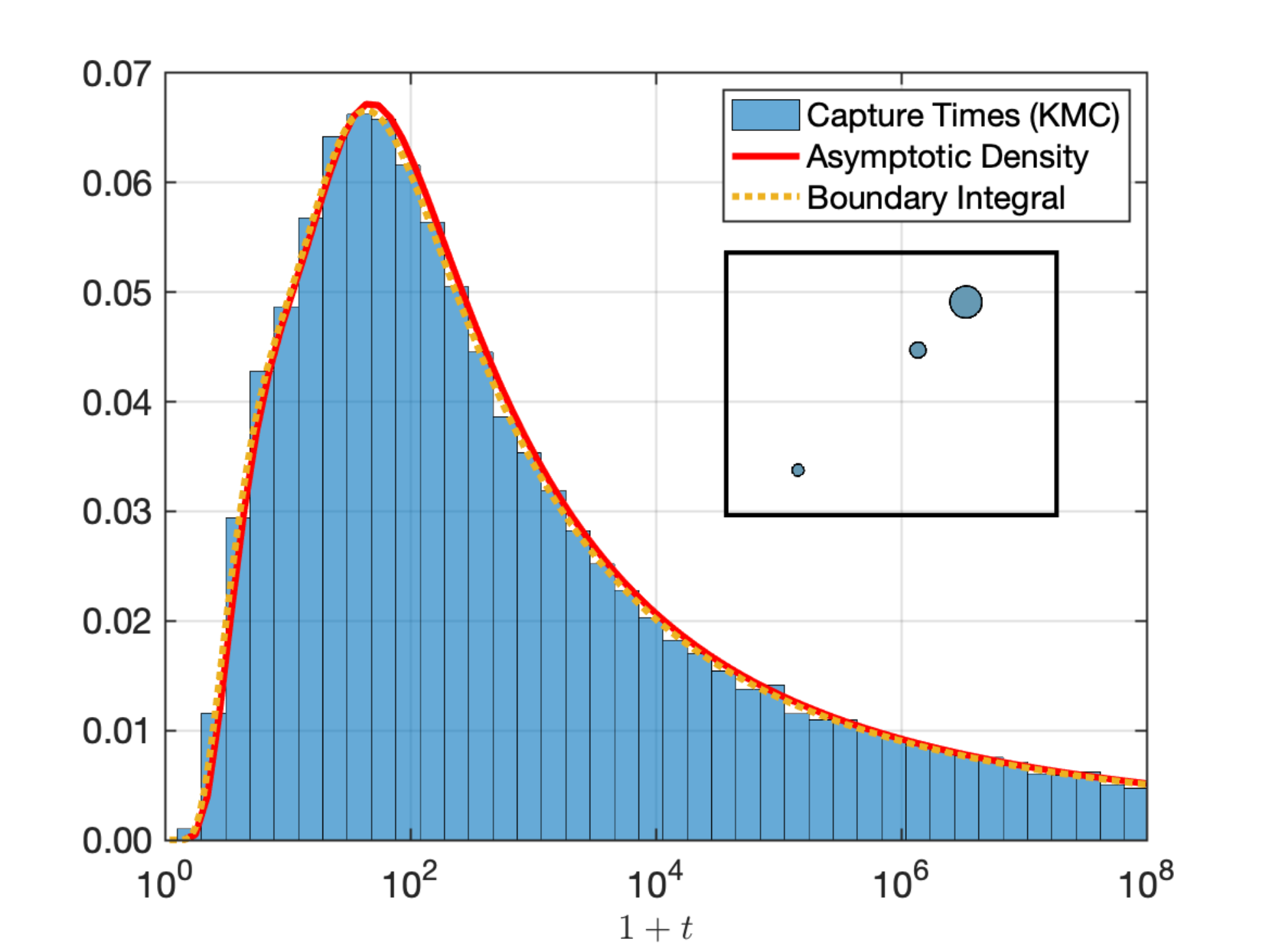}}
   \subfigure[Six target example]{\includegraphics[width=0.325\linewidth]{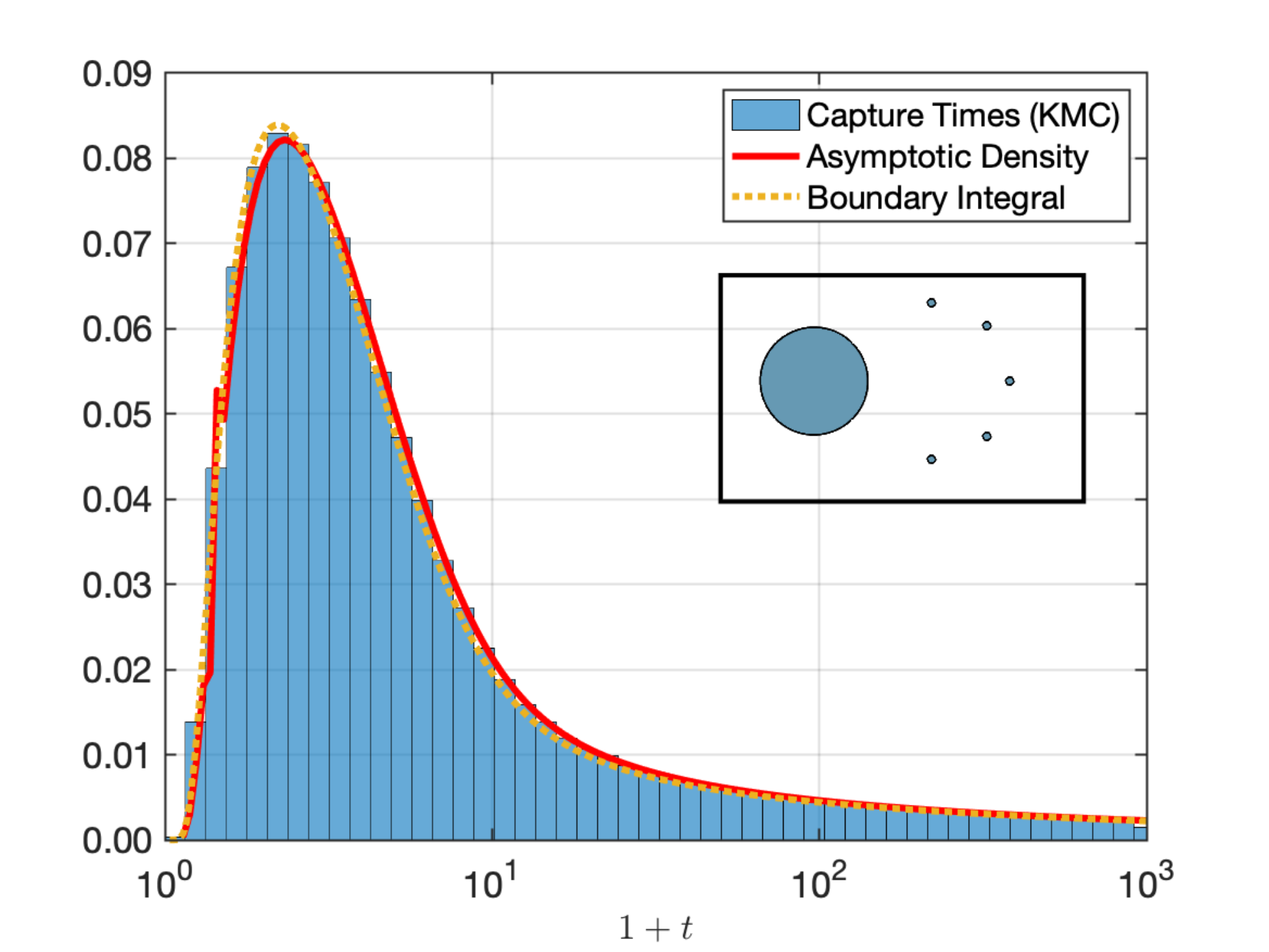}}
      \caption{Relative frequency of arrival times to planar targets by KMC, exact solution (one target case), and hybrid approaches (asymptotic and boundary integral). Schematics of target arrangements shown as insets.\label{fig:PlanarExs} } % include initial position.
% \vspace{-2\baselineskip}% Squeeze space after.
\end{figure}

The results shown in Fig.~\ref{fig:PlanarExs} show good agreement between the three approaches. In Fig.~\ref{fig:PlanarExs}(a) we compare the hybrid-asymptotic, KMC and exact one-pore solutions showing excellent agreement. In the two more challenging examples, we generally see good agreement between the asymptotic and boundary integral approaches. In the more challenging 6 target case, we see in Fig.~\ref{fig:PlanarExs}(c) a slightly diminished agreement is observed near the peak.

\subsection{Homogenization}\label{sec:resultsHomogenization}

In this example, we consider the case of a single target with a mix of absorbing and reflecting portions. This describes the scenario where an impermeable cellular membrane surface is covered in surface receptors.  We determine the full distribution of arrival times using both the KMC method and the hybrid asymptotic-homogenization result \eqref{eqn:homogenization}. In the application of the hybrid approach, we use the analytically determined logarithmic capacitance \eqref{eqn:homogenization}, or equivalently the effective radius, in the single patch result. 

\begin{figure}[htbp]
    \centering
    \subfigure[Ex 1: $M = 12$, $f = 1/3$.]{\includegraphics[width = 0.495\textwidth]{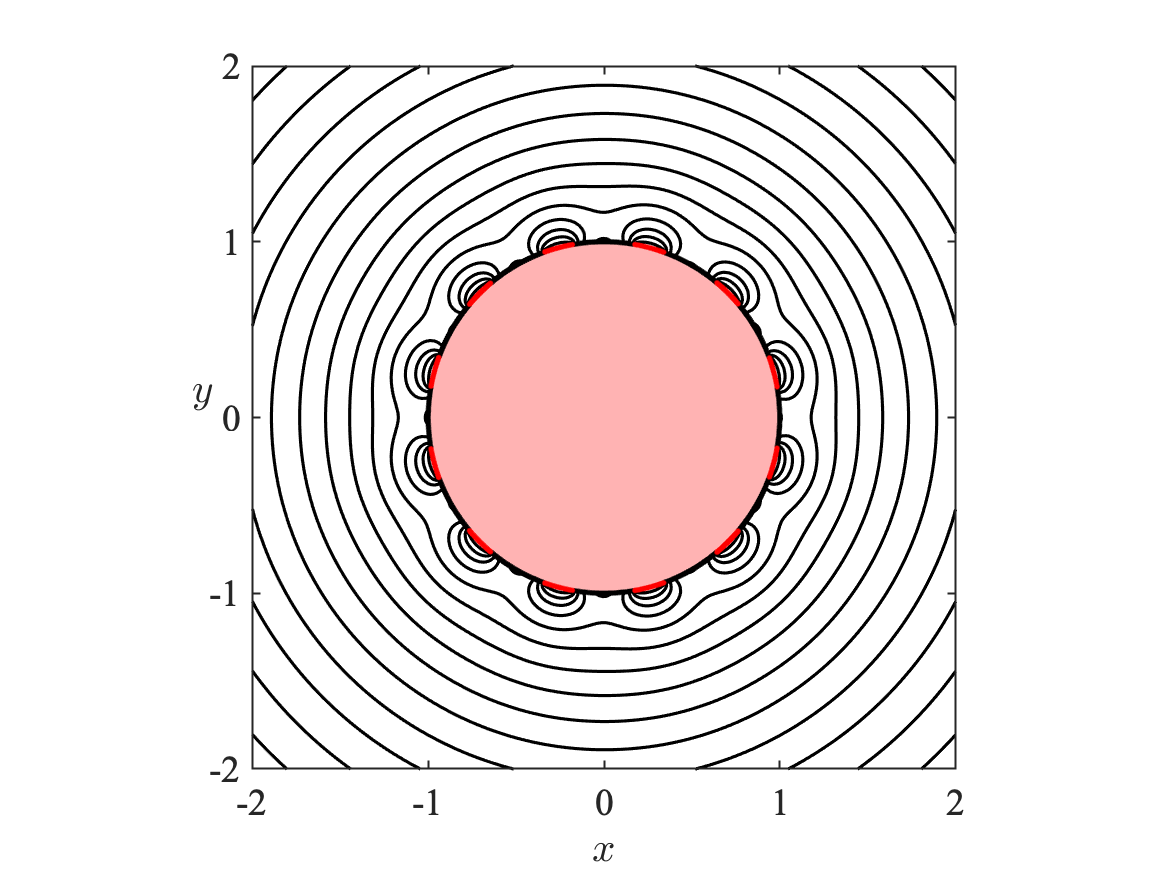}}
    \subfigure[Ex 1: $M = 12$, $f = 1/3$.]{\includegraphics[width = 0.495\textwidth]{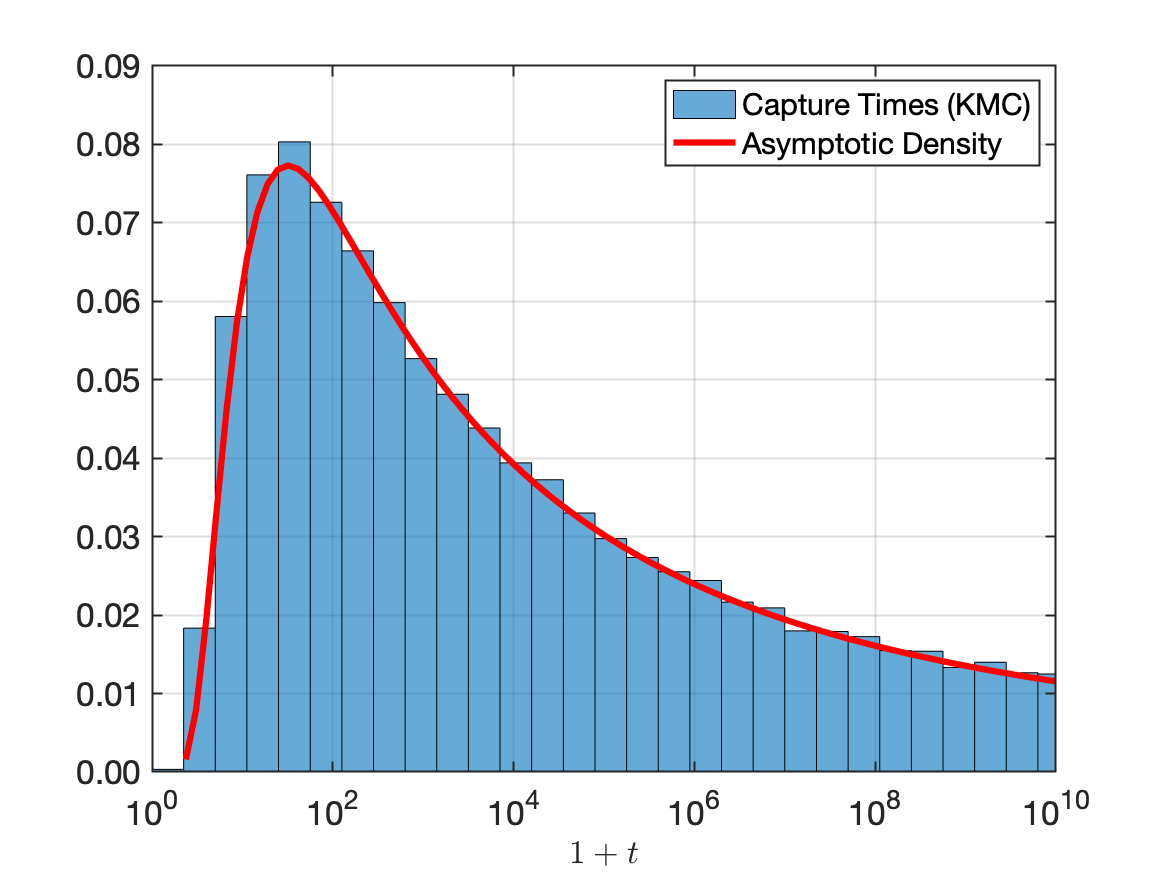}}\\
        \subfigure[Ex 2: $M = 8$, $f = 1/8$.]{\includegraphics[width = 0.495\textwidth]{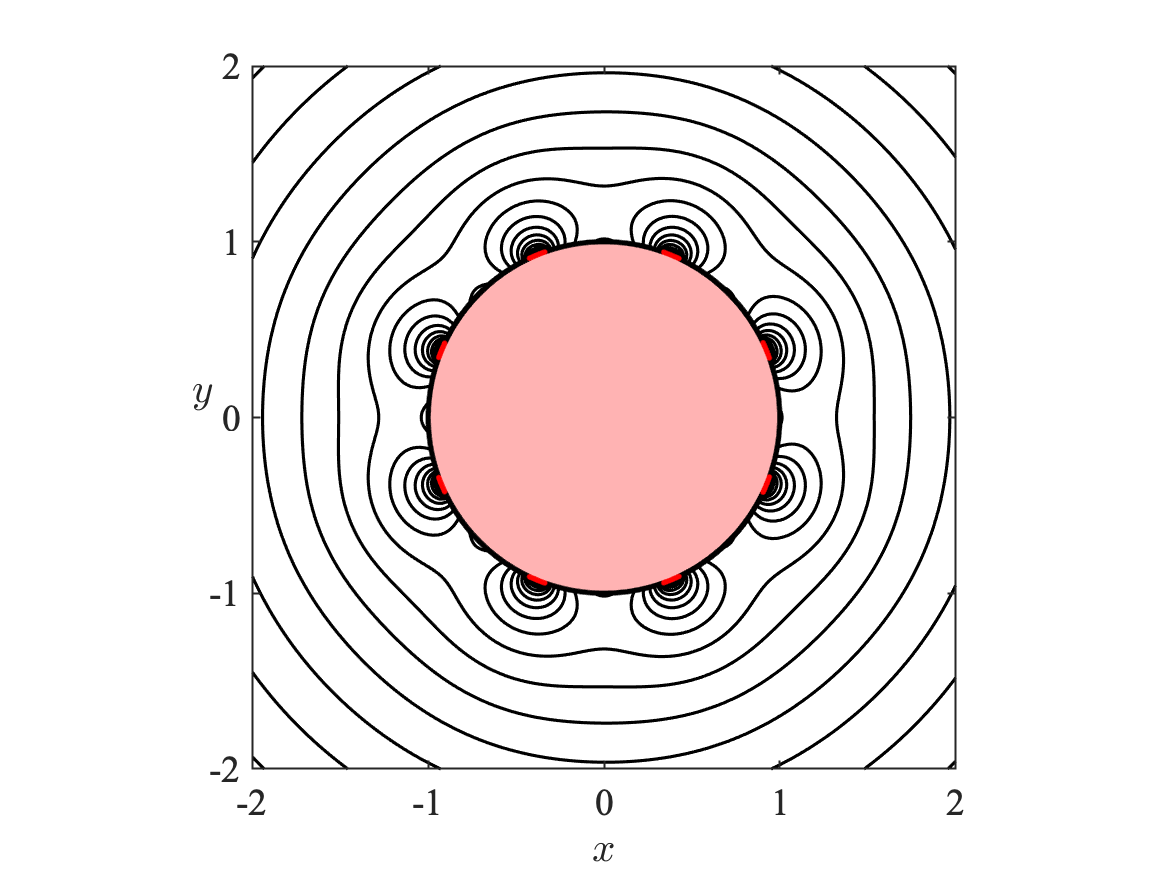}}
    \subfigure[Ex 2: $M = 8$, $f = 1/8$.]{\includegraphics[width = 0.495\textwidth]{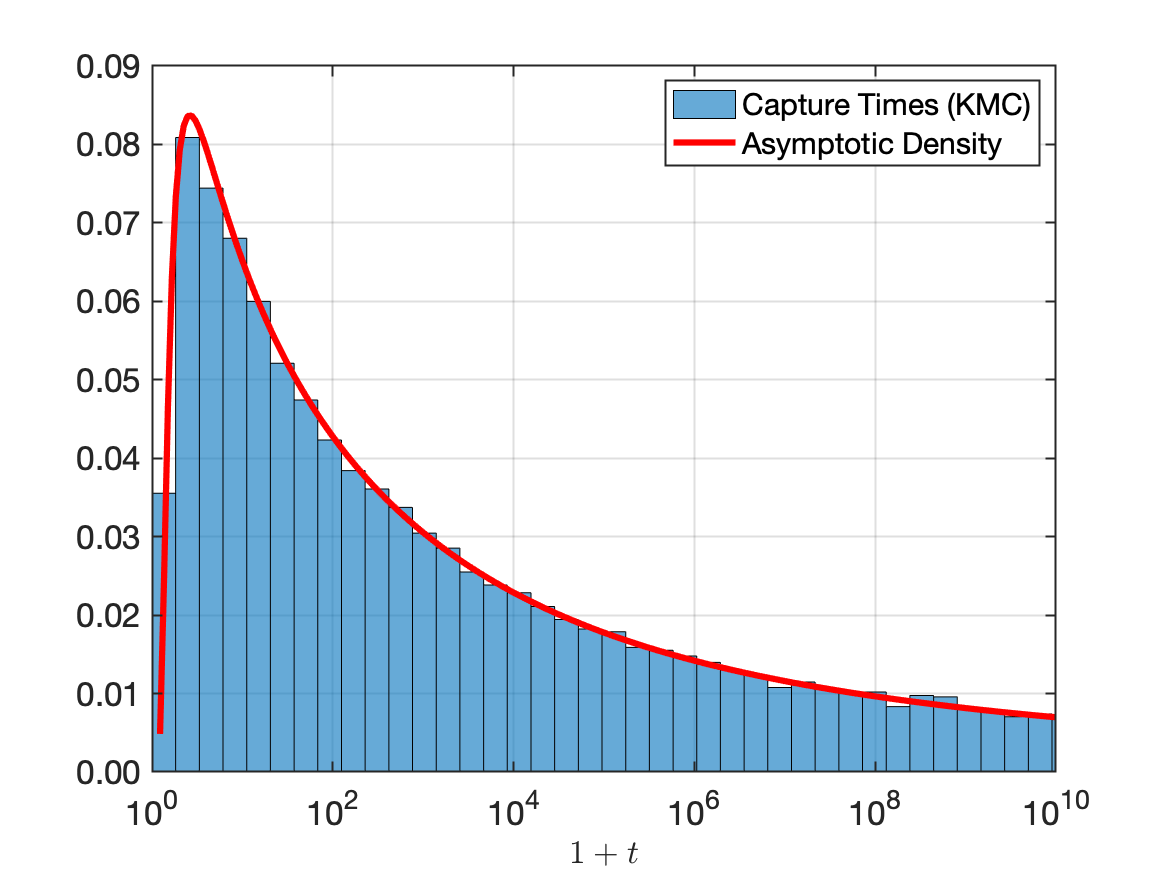}}
    \caption{Results from the KMC and asymptotic method on the homogenization examples \eqref{ex:Homog}. Left panels: target (scaled to unity radius) with layout of absorbers. Contours indicate the numerical solution of \eqref{eqn:LogCap}. Right figures: Agreement between the full arrival time densities obtained from KMC and the hybrid-asymptotic method.}
    \label{fig:HomoG}
\end{figure}

In the two examples shown in Fig.~\ref{fig:HomoG}, we take a circular absorbing target centered at the origin with radius $\eps = 0.05$. The target itself features $M$ equally spaced absorbing windows centered at the roots of unity $(\cos \frac{2\pi k}{M},\sin \frac{2\pi k}{M} )$. The windows occupy a combined fraction $f$ and each has common angular extent $\sigma = \frac{2\pi f}{M}$. In each case we use the homogenized formula \eqref{eqn:homogenization_c} to obtain the logarithmic capacitance $d$ and then apply the result \eqref{2DResults} for $N=1$.  

The relevant parameters obtained for the two examples are
\bsub\label{ex:Homog}
\begin{gather}
\label{ex:Homog_a} \text{Ex 1:} \quad M = 12, \quad f = \frac{1}{3}, \quad d = 0.8978, \quad \bx_0 = [5,0]; \\[5pt]
\label{ex:Homog_b} \text{Ex 2:} \quad M = 8,\quad f =  \frac{1}{8}, \quad  d = 0.6657, \quad \bx_0 = [2,0].
\end{gather}
\esub
In Fig.~\ref{fig:HomoG}(b),(d) we show good agreement between the two methods while in Fig.~\ref{fig:HomoG}(a),(c) we display a visualization of the solution to the capacitance problem \eqref{eqn:LogCap}. We remark that the homogenization effect can be seen through the gradual radial symmetrization of the contours. 

\subsection{Clustering}\label{sec:resultsClustering}

In this section, we use two examples of clustered target configurations to compare results from the KMC method with both the asymptotic and BIM approaches. In addition, we show the effectiveness of homogenization where the clustered target configuration is replaced by a single circular target of appropriately chosen radius. The specific parameters are given by
\bsub\label{ex:cluster}
\begin{align}
\text{Ex 1:} \quad \bx_{k} &=  \frac{3}{2}\left(\cos \frac{2\pi k}{5},\sin \frac{2\pi
  k}{5}\right), \quad k = 1,\ldots,5, \quad \bx_6 = (0,0), \qquad  r_{1,\ldots,6} = 0.1, \quad \bx_0 = (10,0);  \\[5pt]
\text{Ex 2:} \quad \bx_{k} &=  \frac{3}{2}\left(\cos \frac{2\pi k}{8},\sin \frac{2\pi
  k}{8}\right), \quad k = 1,\ldots,8, \qquad r_{1,\ldots,8} = 0.05, \quad \bx_0 = (5,0).
\end{align}
\esub

\begin{figure}[htbp]
\centering
\subfigure[Numerical solution of \eqref{eqn:LogCap} yielding logarithmic capacitance $d = 1.2295$.]{\includegraphics[width = 0.23\textwidth]{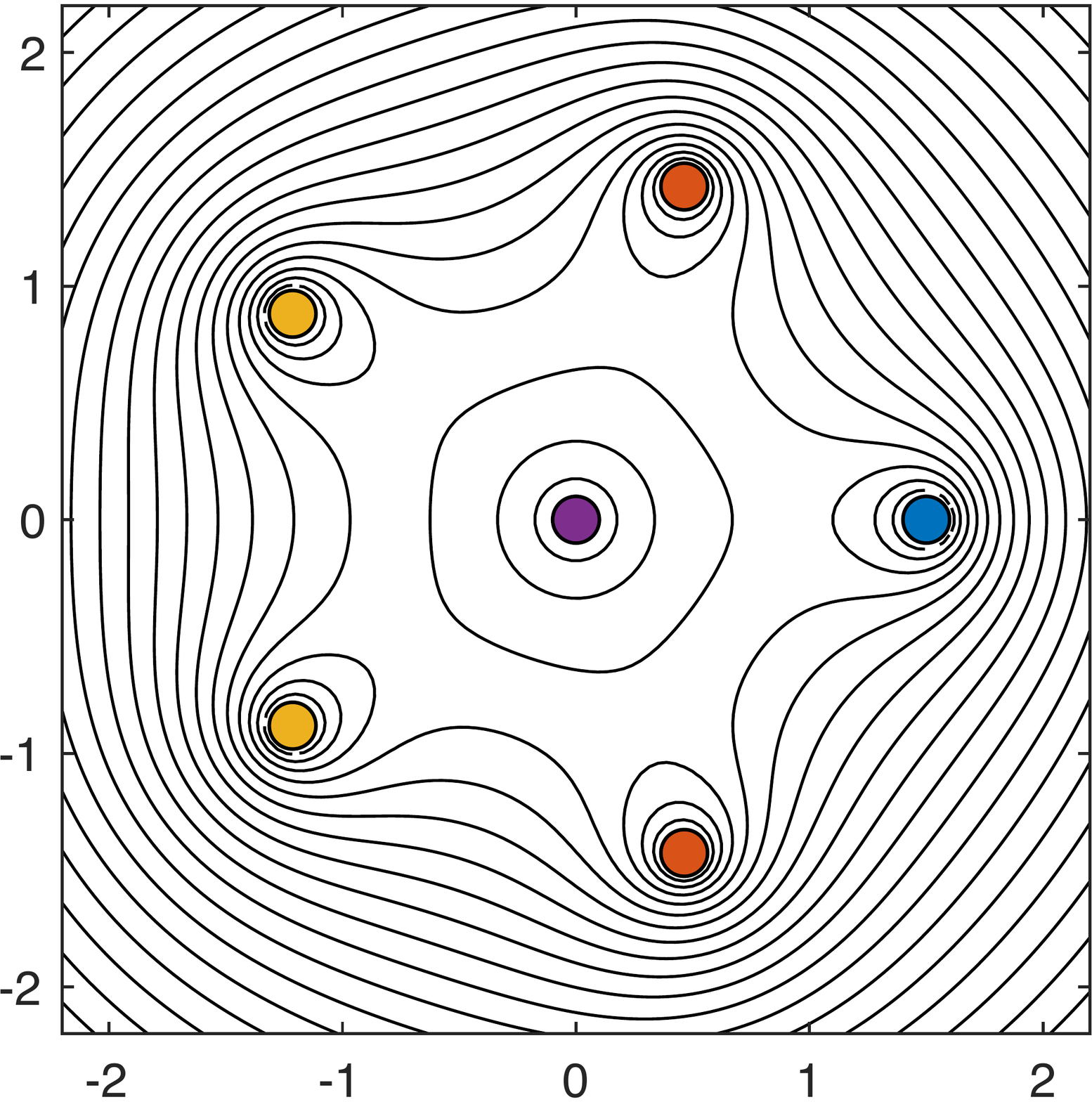}}\quad
\subfigure[Normalized target fluxes.]{\includegraphics[width = 0.32\textwidth]{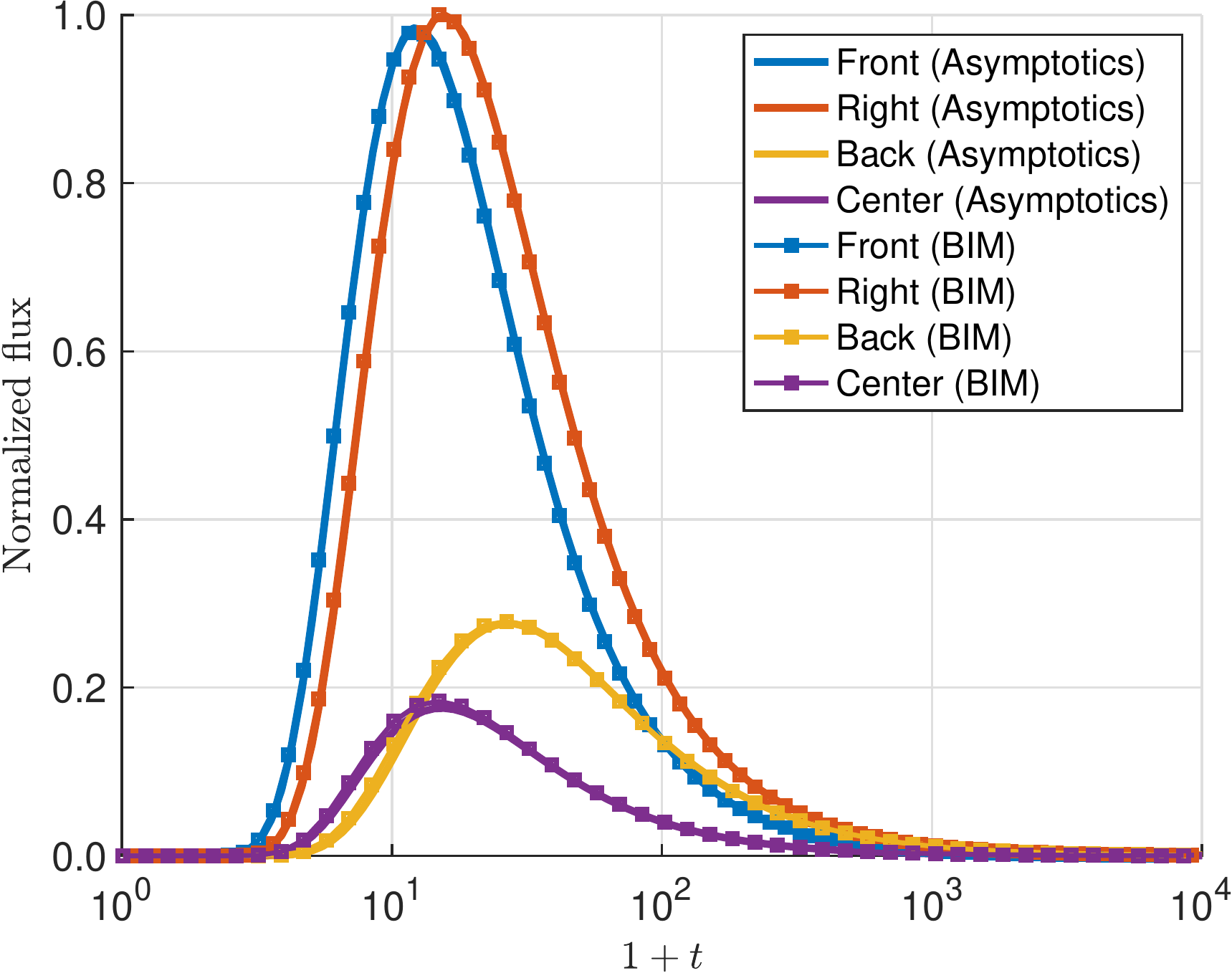}}\quad
\subfigure[Agreement of capture times (from KMC) with asymptotic and homogenized densities.]{\includegraphics[width = 0.35\textwidth]{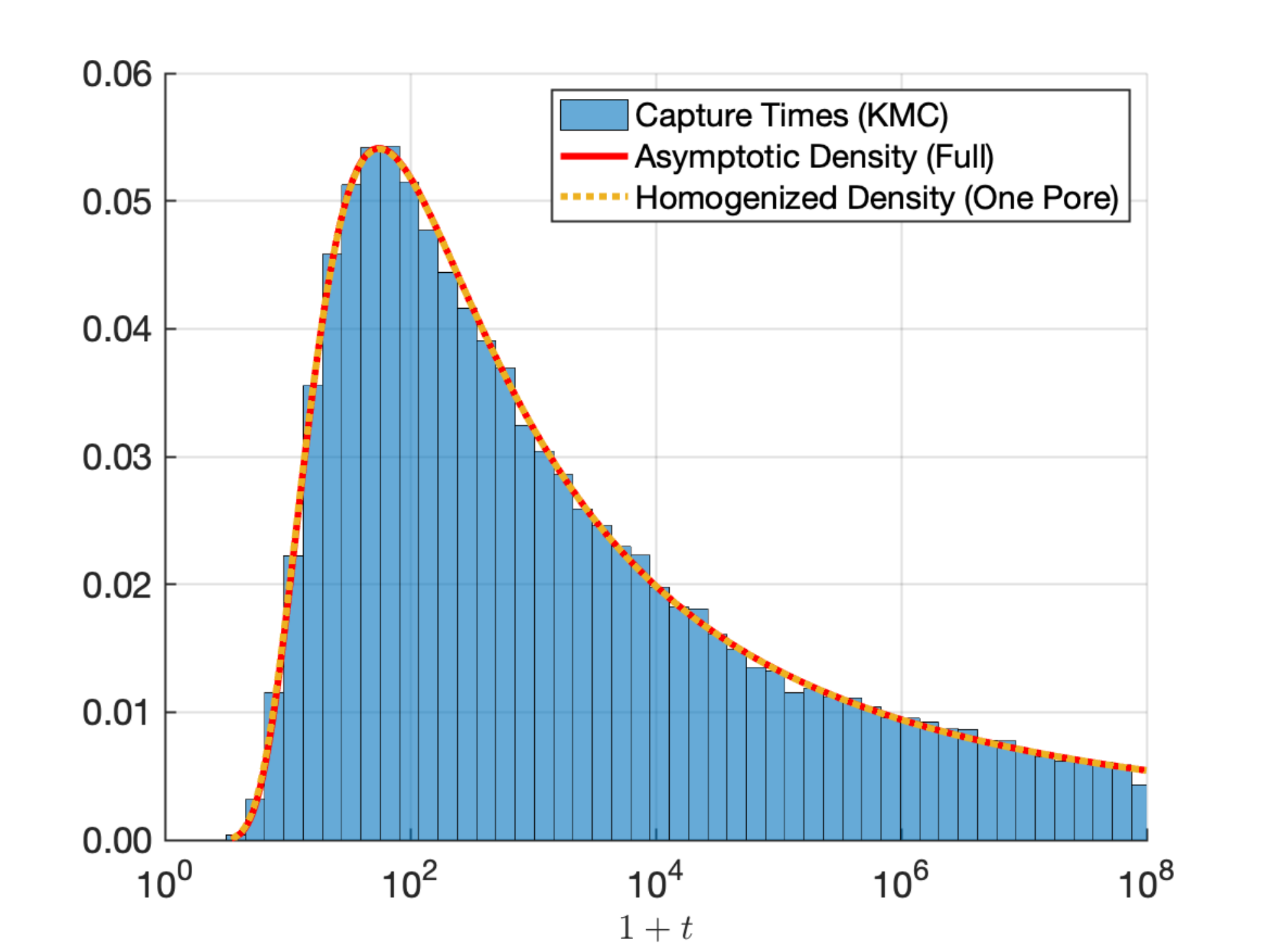}}\\
%%%%%%%%%%%%%%%%%%%%%%%%%%%%%%%%%%%%%%%%%%%%%%%%%%%%%%%%%%%%%%%%%%%%%%%%%%%%%%%%%%%%
\subfigure[Numerical solution of \eqref{eqn:LogCap} yielding logarithmic capacitance $d = 1.2737$.]{\includegraphics[width = 0.24\textwidth]{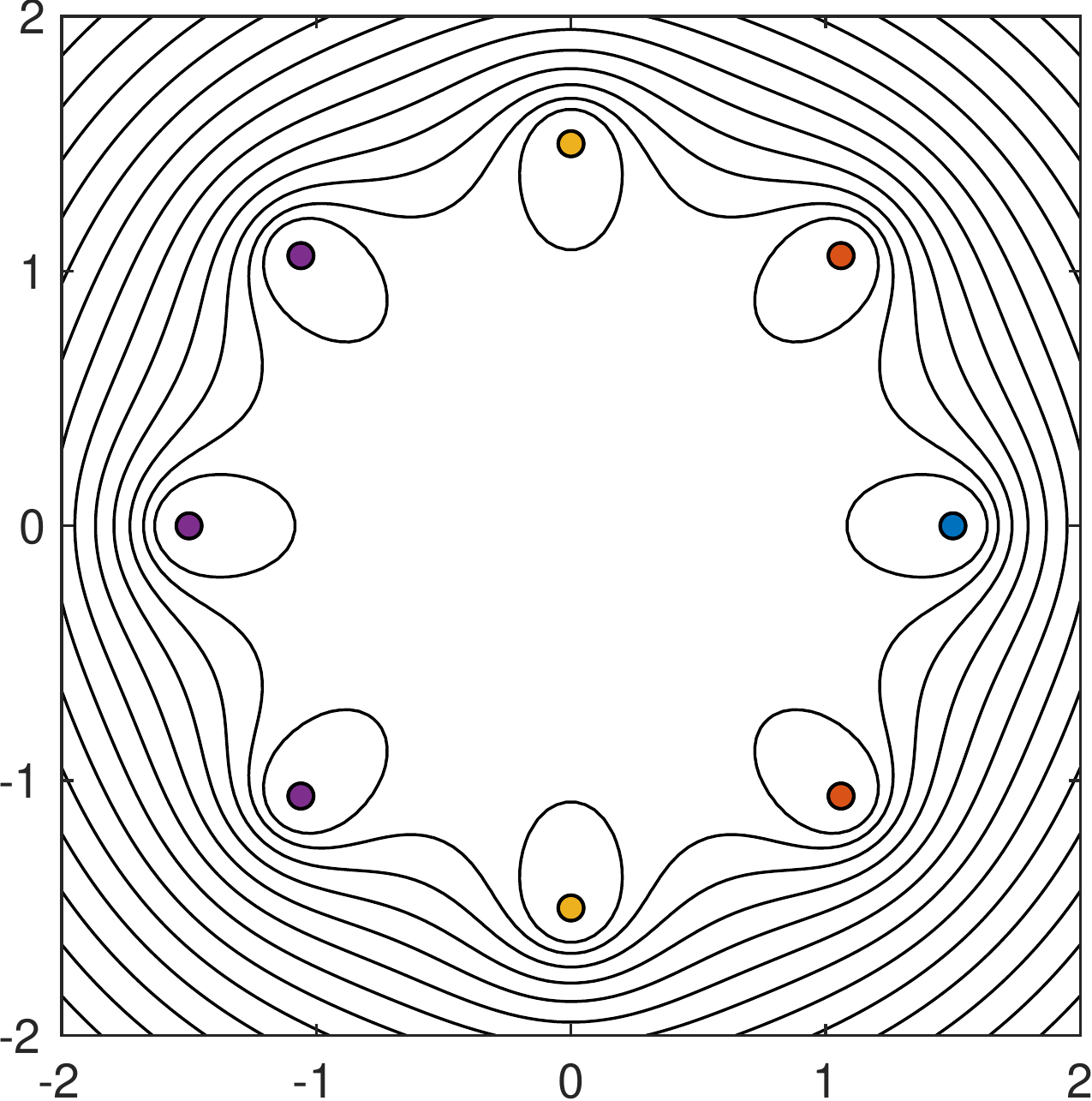}} \quad  
\subfigure[Normalized target fluxes.]{\includegraphics[width = 0.32\textwidth]{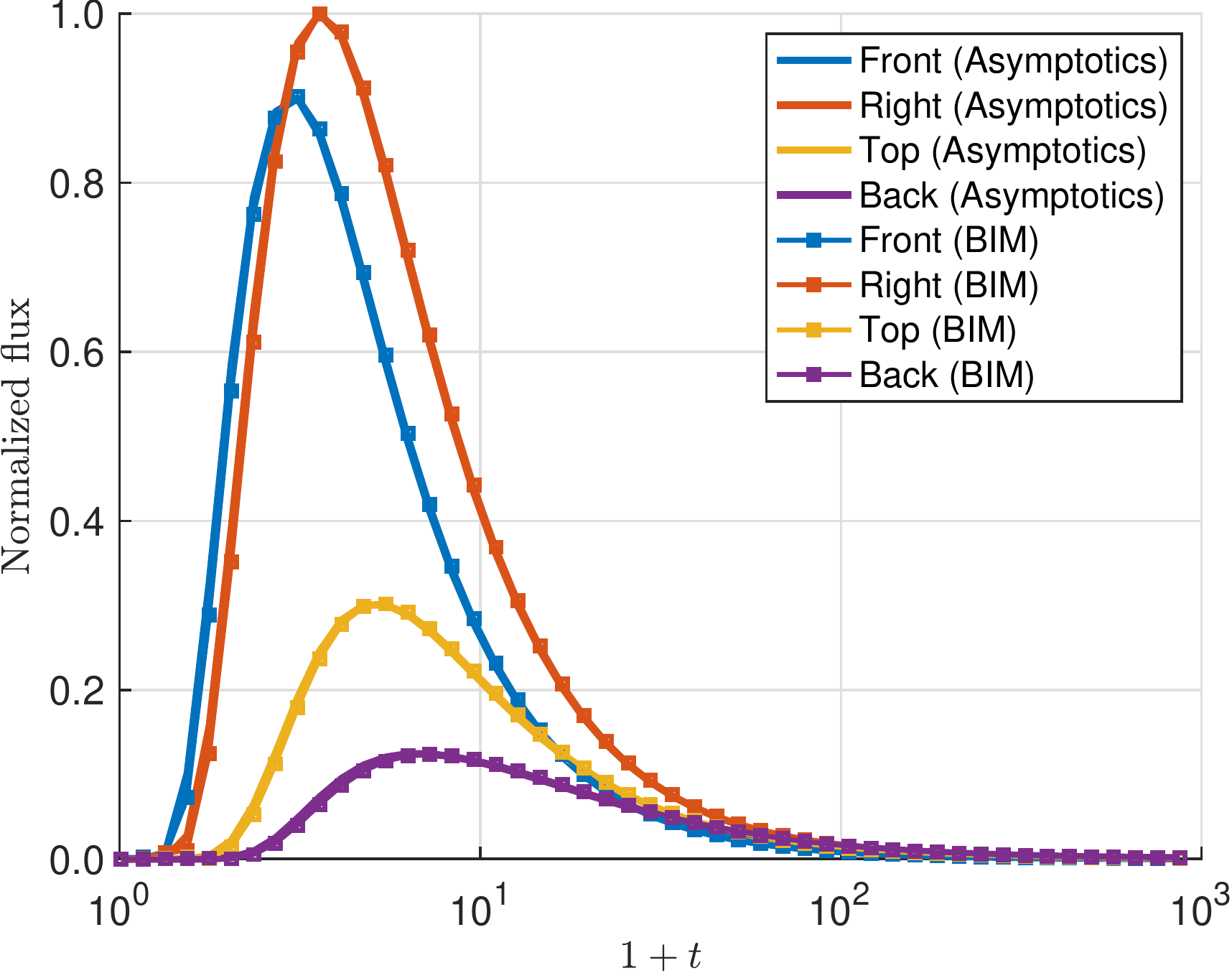}}\quad
\subfigure[Capture times (from KMC) with asymptotic and homogenized densities.]{\includegraphics[width = 0.35\textwidth]{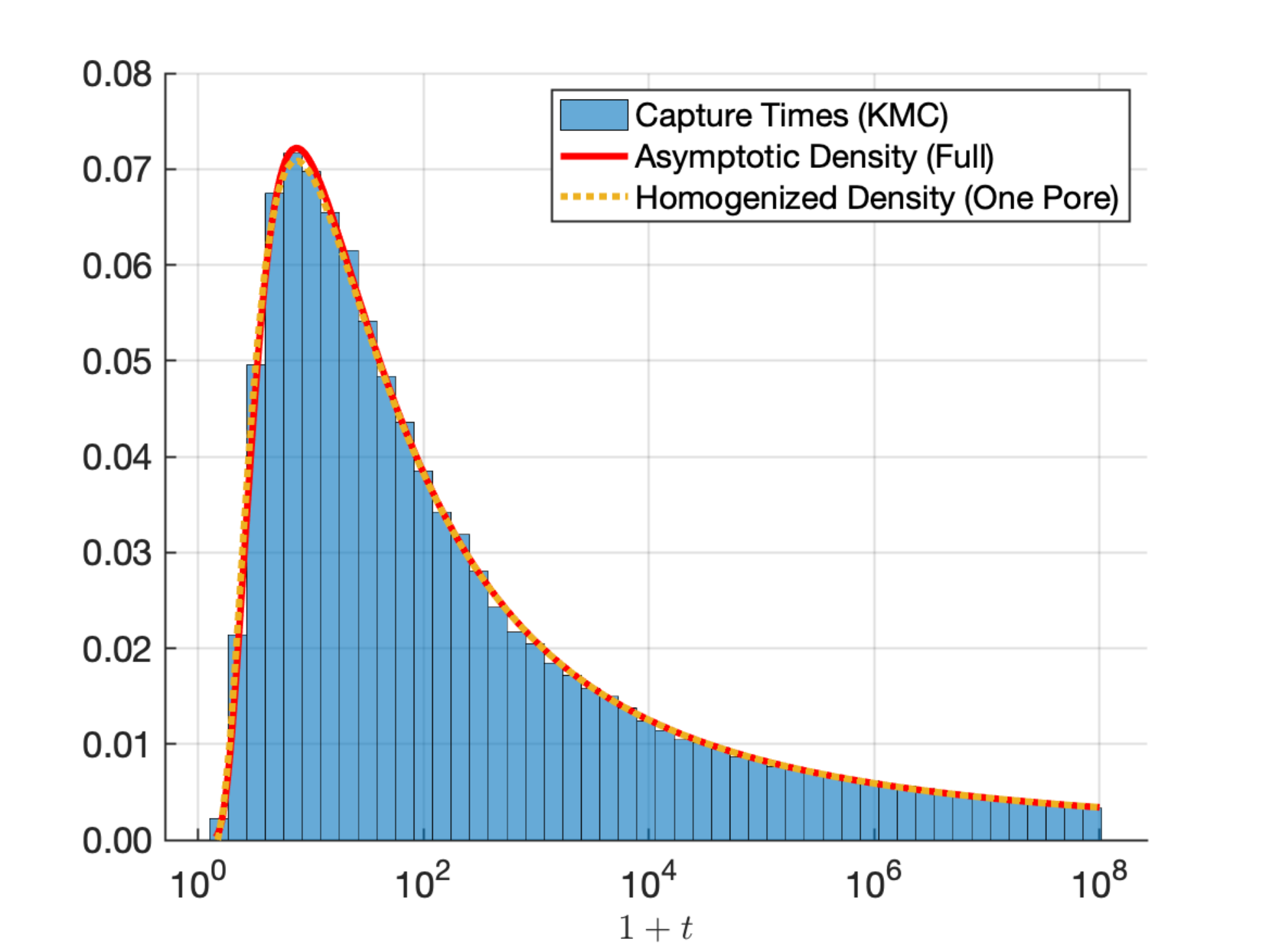}}
\caption{Comparison of the full asymptotic density with those arising from the KMC method, boundary integral method (BIM) and homogenization (replacing the target cluster with a single absorbing target). }
\label{fig:Clustering}
\end{figure}
The numerical method described by equation \eqref{eqn:complexSeries} allows for the computation of the logarithmic capacitance for each of the absorbing sets (cf.~Fig.~\ref{fig:Clustering}(a,d)). For the parameters specified in \eqref{ex:cluster}, we determine $d = 1.2295$ (\text{Ex 1}) and $d= 1.2737$ (\text{Ex 2}). The logarithmic capacitance can be interpreted as the effective radius of a single target that reflects the capture potential of the cluster. As with the homogenization example in Sec \ref{sec:resultsHomogenization}, we observe that replacing complex configurations with a single target of appropriately chosen radius produces a very accurate representation of the full arrival time distribution. However, the single target representation does reduce certain direction information encoded in the distribution of arrivals over the targets in the cluster. The hybrid-asymptotic method can rapidly determine the fluxes by numerical inverse Laplace transform of $\hat{\mathcal{J}}_j =  2\pi D \nu_j S_j$  where each $S_j$ is determined by \eqref{sysmain}. In Fig.~\ref{fig:Clustering}(b,e) we show normalized fluxes into four sets of traps that are symmetrically arranged with respect to the initial location. The traps aligned towards the initial data accrue most of the inbound flux while the peaks, representing most likely arrival times at a particular target, are ordered by their distance to the initial location. The distribution of fluxes over the targets encode directional information that can infer the source location \cite{Miles2020}.

\subsection{Splitting Probabilities}\label{sec:resultsSplitting}

In this section, we demonstrate the convergence of the dynamic fluxes to the static splitting probabilities $\{\phi_k(\bx)\}_{k=1}^N$, where $\phi_k(\bx)$ denotes the probability that a diffusing particle originally at $\bx$ reaches the $k^{th}$ target before any others. For exposition purposes, we focus on the scenario of completely absorbing targets. These probabilities satisfy the exterior Laplace problem
\bsub\label{eq:splits}
\begin{gather}
\label{eq:splits_a} \Delta \phi_k = 0, \quad \bx\in\mathbb{R}^2\setminus\Omega; \qquad \phi_k(\bx) = \delta_{jk},\quad \bx \in\partial\mathcal{A}_j, \qquad j= 1,\ldots, N;\\[5pt]
\label{eq:splits_b} \phi_{k}(\bx) \quad \text{finite as} \quad |\bx| \to \infty,
\end{gather}
\esub
where $\delta_{jk}$ is the Kronecker delta. The asymptotic solution of \eqref{eq:splits} as $\eps\to0$ is developed along similar lines to Sec.~\ref{sec:asy} (also see \cite[Sec 5]{Kurella2015}). Accordingly, we present the solution as $\eps\to0$ directly as
\bsub\label{eqn:splitsol}
\begin{equation}\label{eqn:splitsol1}
\phi_k(\bx) \sim 2\pi \sum_{j=1}^{N} \nu_j A_j G_0(\bx;\bx_j) + \bar{\phi}, \qquad G_0(\bx;\by) = \frac{1}{2\pi}\log|\bx-\by|,
\end{equation}\label{eqn:splitsol2}
where the $(N+1)$ constants $(A_1,\ldots,A_N, \bar{\phi})$ are determined from the linear system,
\begin{equation}
\sum_{j=1}^{N} \nu_j A_j = 0; \qquad -A_j + 2\pi \sum_{\substack{i=1\\ i\neq j}}^N \nu_i A_i G_0(\bx_j;\bx_i)  + \bar{\phi} = \delta_{jk}, \qquad j = 1,\ldots,N.
\end{equation}
\esub
The gauge functions $\nu_{j} = -1/\log\eps d_j$ are defined in \eqref{eq:guage}. We remark that since capture is guaranteed for planar Brownian motion, we have that $\sum_{k=1}^N \phi_k(\bx) = 1$.

\begin{figure}[htbp]
\centering
\subfigure[Schematic of target configuration.]{\includegraphics[width=0.45\textwidth]{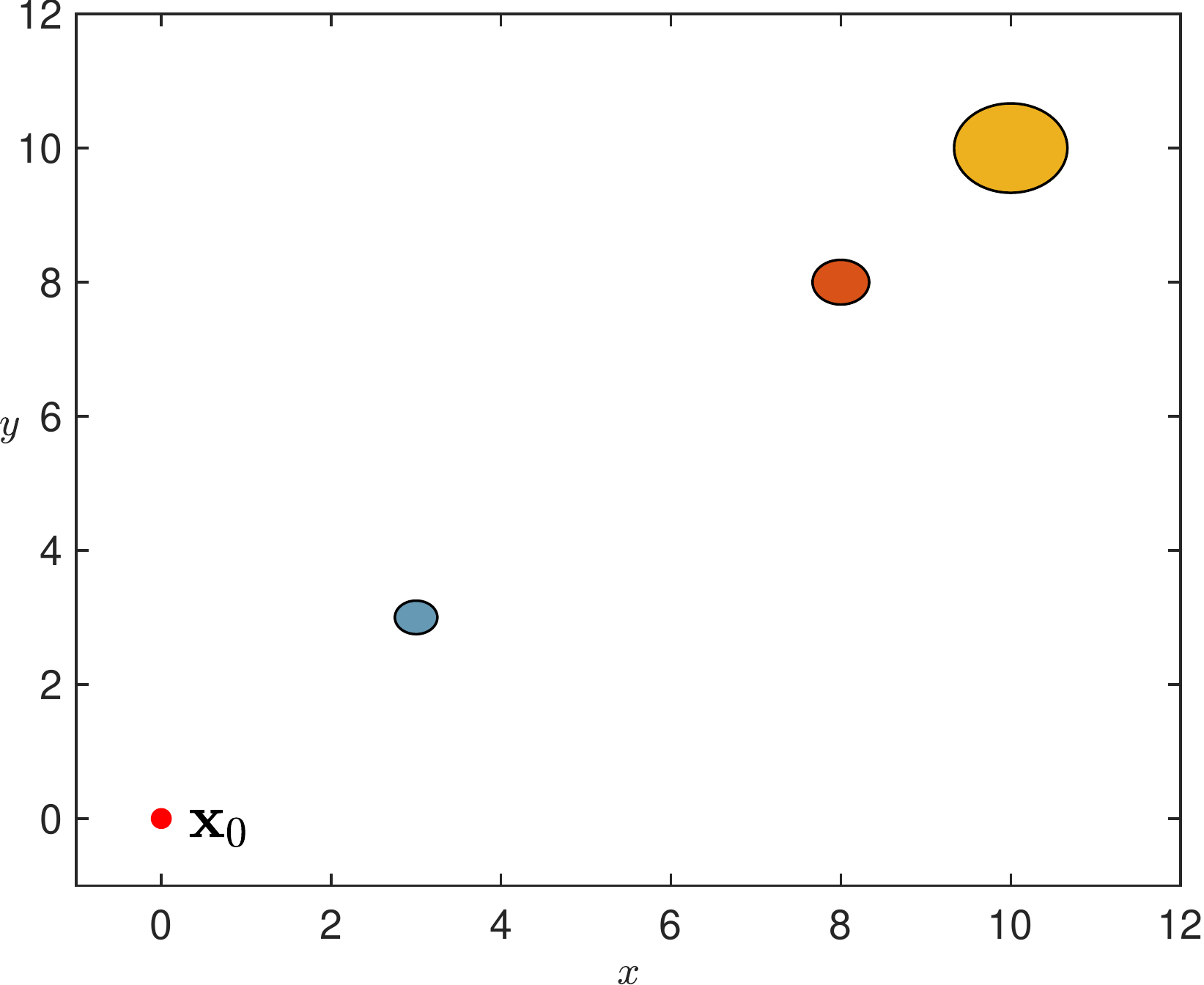}}\qquad
\subfigure[Splitting probabilities and dynamic fluxes.]{\includegraphics[width=0.45\textwidth]{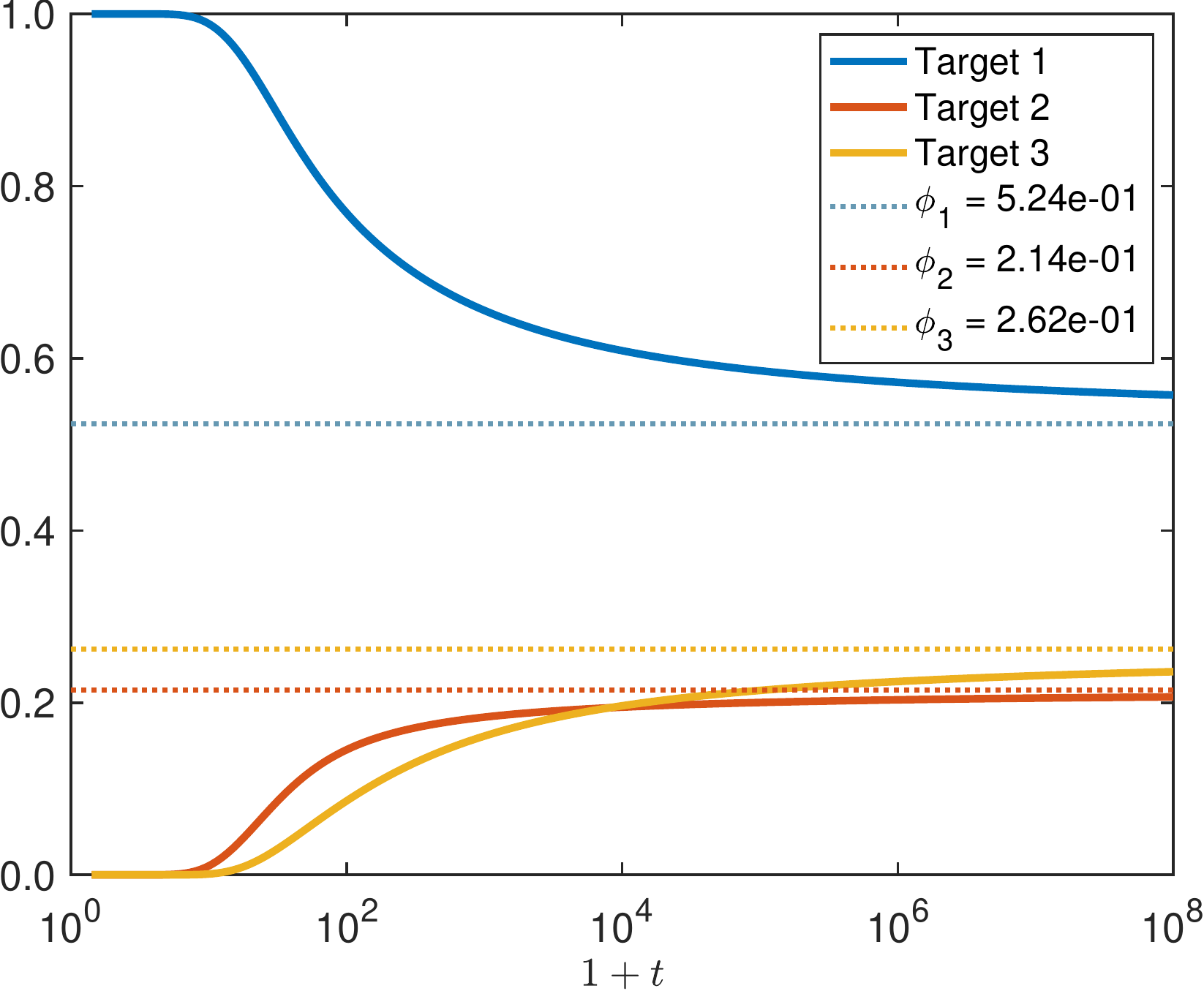}}
\caption{Left: A schematic of the three target configuration from Sec.~\ref{sec:resultsPlanar}. Right: The fractional fluxes $q_k(t) = \mathcal{J}_k(t)/\sum_{j=1}^N \mathcal{J}_j(t)$ into each target together with the steady state splitting probabilities $\phi_k$ obtained from \eqref{eqn:splitsol}. \label{fig:splits}}
\end{figure}
To demonstrate this theory, we compare these static splitting probabilities with the time-dependent fractional fluxes $q_k(t) = \mathcal{J}_k(t)/\sum_{j=1}^N \mathcal{J}_j(t)$ into each target (see Fig.~\ref{fig:splits}) obtained from the hybrid-asymptotic method (Sec.~\ref{sec:asy}). We draw attention to two important conclusions from this example. First, the dynamic fluxes converge to the static splitting probabilities on a very long timescale. For many physically or biologically relevant timescales, this questions the usefulness of using splitting probabilities for inference \cite{Miles2020} purposes. Second, the ordering of the relative fluxes into each target changes over the displayed time interval. Specifically, at short times ($t\lessapprox10^1$), target 1 captures almost all the flux and indeed is the most significant absorber over the entire timeline. Target 1 is the smallest target but closest to the initial location demonstrating that this distance is a significant indicator of capture ability. Later in the timeline, we see that targets 2 and 3 interchange their prominence in capture fluxes (around $t\approx 10^4$) implying that proximity to the initial location promotes faster capture at short times while at larger times, target size can play a more significant role. Importantly, none of these subtleties are apparent from the static splitting probabilities highlighting the necessity of obtaining full time dependent statistics.

\section{Discussion}
% method, applications, and further work.
In this work we have demonstrated several methodologies for obtaining arrival time statistics of diffusing particles to complex sets of absorbing targets in planar regions. The Laplace transform approach seeks to solve an equation of modified Helmholtz type by either an asymptotic expansion for well separated targets or a boundary integral representation. In both cases the target geometries can be very general, but for the boundary integral method, it is presently limited to purely absorbing targets. The inverse transform is obtained by quadrature of the Bromwich integral. To complement these approaches, we developed a particle based kinetic Monte-Carlo (KMC) method that can resolve the arrival distribution for very general configurations of targets and boundary sets. These methods are rapid, accurate and easy to implement. The hybrid asymptotic method is particularly suited to the scenario of well-separated targets while the KMC method is applicable to general geometric scenarios and varied boundary conditions. 

\begin{figure}
\centering
\includegraphics[width = 0.45\textwidth]{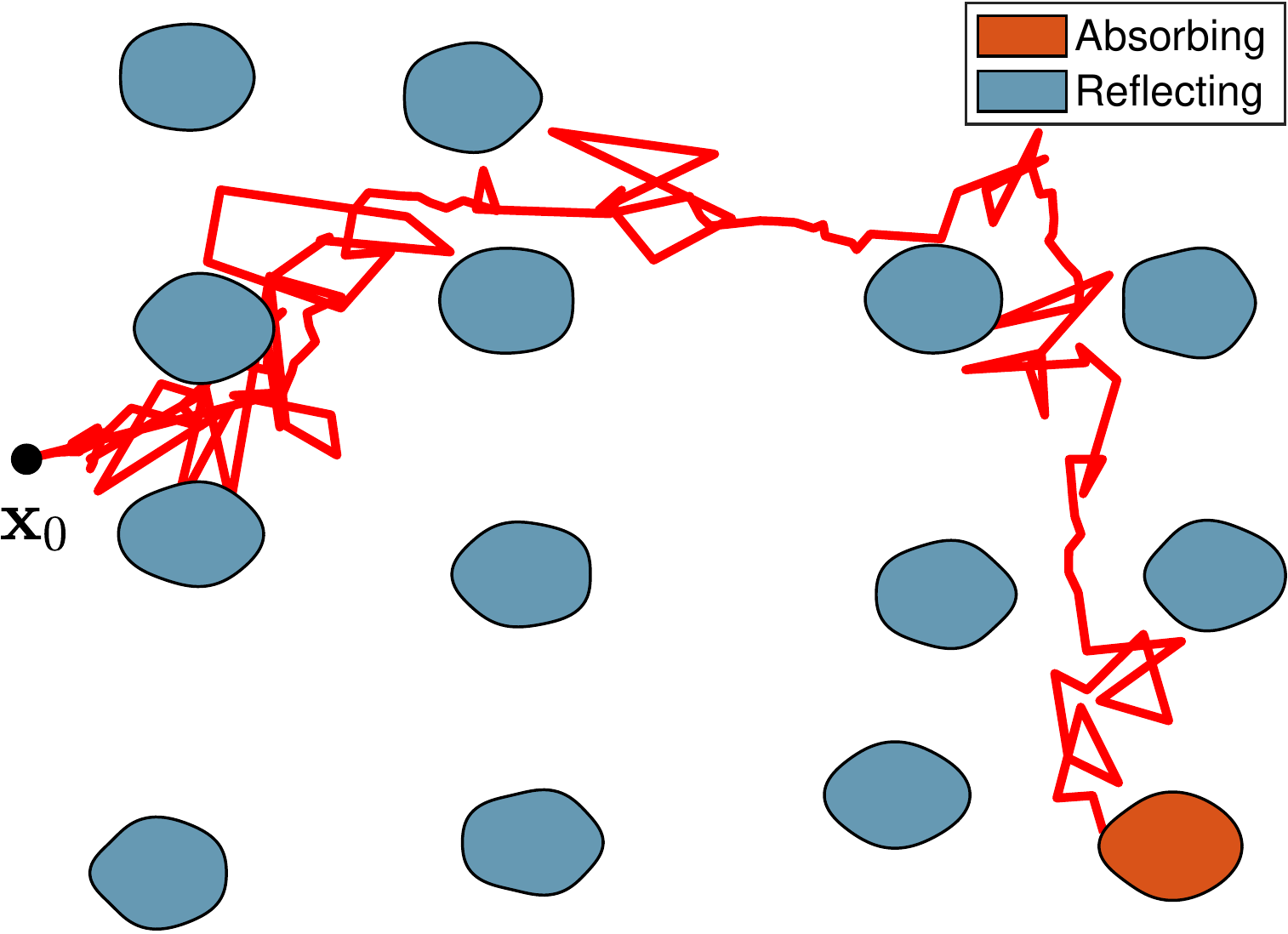}
\caption{A tricky first passage time distribution to evaluate using either KMC or hybrid-asymptotic methods.\label{fig:crazypath}}
\end{figure}

There are a few conclusions that emerge from our study. Homogenization is a powerful technique that accurately reproduces the first passage time distributions of complex target sets by replacing them with a single circular target of appropriately chosen radius. However, homogenization brings limitations with it, particularly as it coarse grains the spatial distribution of arrivals across the targets. The relative fraction of particles that arrive across a distribution of targets has directional information that can be used to infer source location \cite{Miles2020}. Additionally, the dynamics available from the full distribution of arrival statistics reveals the limitations of using static information, (e.g.~splitting probabilities) which are only representative of very long time behavior. Our example of dynamic splitting probabilities (Sec.~\ref{sec:resultsSplitting}) suggests that there are several timescales over which arrival statistics can be important and that the relevant physical or biological timescale must be considered when inferring source location from arrival information.

While the methods developed here are quite general, there are certain scenarios where they fail and new approaches are needed. For example, in the scenario where there are numerous targets but only a few are reactive (Fig.~\ref{fig:crazypath}), particles must navigate a torturous route through inert targets to reach the destination. The combination of targets with either purely Neumann or Dirichlet boundary conditions hampers all the approaches developed here. This scenario is particularly challenging for the KMC method due to the fact that on a reflecting target surface, the particle will tend to perform a surface diffusion characterized by many small jumps. The very long wait time for a large jump necessary to leave the target vicinity means that the convergence rate is greatly reduced. To address these shortcoming, we plan extensions to the boundary integral formulation. 

\appendix

\section{Two dimensional problem arrival problem. Arrival time and angle distribution}\label{sec:PlanarDist}
For a particle with diffusivity $D$ originally at $\bx\in(b, 0)$, the occupation density in $r=|\bx|>a$ satisfies
\bsub\label{eqn:appstart}
\begin{gather}
    p_t = D\left(p_{rr} + \frac{1}{r} p_r + \frac{1}{r^2}
    p_{\theta\theta}\right) ,\qquad  r>a, \quad \theta \in(-\pi,\pi), \quad t>0;  \\[5pt]
    p(a,\theta,t) = 0,\qquad \theta\in(-\pi,\pi), \quad t>0;\\[5pt]
    p(r,\theta,0) = \frac{1}{r}\delta(r-b)\delta(\theta), \qquad r>a, \quad \theta\in(-\pi,\pi).
\end{gather}
\esub
Our goal is to determine closed form expressions for the quantities
\bsub
\begin{gather}
P(t) = \int_{\theta=-\pi}^{\pi} \int_{r=a}^{\infty} p(r,\theta,t) \, r drd\theta, \qquad [\text{Survival probability}] \\
C(t) =1-\int_{\theta=-\pi}^{\pi} \int_{r=a}^{\infty} p(r,\theta,t) \, r drd\theta, \qquad [\text{Capture probability}] \\
 S(t) = - P'(t), \qquad [\text{Arrival time distribution}] 
\end{gather}
\esub
Applying the divergence theorem to $S(t) = -P'(t)$, we see that
\begin{equation}
S(t) = - \int_{\theta=-\pi}^{\pi} \int_{r=a}^{\infty} p_t(r,\theta,t) \, r drd\theta = - D \int_{\theta=-\pi}^{\pi} \int_{r=a}^{\infty}( r p_r )_r \, dr d\theta = D a \int_{\theta=-\pi}^{\pi}  p_r(a,\theta,t) \, d\theta.
\end{equation}
\bsub\label{eqn:rescaled}
To obtain the flux $p_r(a,\theta,t)$, we non-dimensionalize by introducing variables
\begin{equation}\label{eqn:rescaled_vars} 
p(r,t) =  \frac{1}{a^2}\, \tilde{p}(\tilde{r},\theta,\tilde{t}),\qquad  \tilde{r} = \frac{1}{a} r, \qquad \tilde{t} =  \frac{D}{a^2} t, \qquad R = \frac{b}{a}.
\end{equation}
so that $S(t) = \frac{D}{a^2}\int_{\theta=-\pi}^{\pi}
\tilde{p}_{\tilde{r}}(1,\theta,\tilde{t})\,d\theta$.
Under the change of variables \eqref{eqn:rescaled_vars}, \eqref{eqn:appstart} becomes
\begin{gather}
\label{eqn:rescaled_a}    \tilde{p}_{\tilde{t}}  = \left( \tilde{p}_{
  \tilde{r}\tilde{r} } + \frac{1}{ \tilde{r} }  \tilde{p}_{\tilde{r}} +
\frac{1}{ \tilde{r}^2} \tilde{p}_{\theta\theta} \right) ,\qquad  \tilde{r}> 1, \quad \theta \in(-\pi,\pi), \quad \tilde{t}>0;  \\[5pt]
\label{eqn:rescaled_b}    \tilde{p}(1,\theta, \tilde{t}) = 0,\qquad  \theta\in(-\pi,\pi), \quad \tilde{t}>0;\\[5pt]
\label{eqn:rescaled_c}    \tilde{p}(\tilde{r},\theta,0) =  \frac{1}{\tilde{r}}\delta(\tilde{r}- R )\delta(\theta), \qquad \tilde{r}>1, \quad \theta\in(-\pi,\pi).
\end{gather}
\esub
After dropping the tildes, we solve for the dimensionless occupation
density $p(r, \theta, t)$ by transforming to Laplace space
$\hat{p}(r,\theta,s)= \int_{t=0}^{\infty}p(r,\theta,t) \,e^{-st}dt$, to see that \eqref{eqn:rescaled_a} satisfies the PDE
\begin{equation}\label{eq:disc2}
    \hat{p}_{rr}+\frac{1}{r}\hat{p}_{r}+\frac{1}{r^2}\hat{p}_{\theta\theta}-sp=-\frac{1}{r}\delta(r-R)\delta(\theta),
    \qquad r>1, \quad \theta \in(-\pi,\pi), \quad s\in\CC.
\end{equation}
The separable solution that is continuous, bounded and satisfies $\hat{p}(1,\theta,s)=0$ is
\begin{equation}
\hat{p}(r,\theta,s)= \left\{  \begin{array}{lc}
\ds\sum_{n=0}^{\infty} A_n\left[ I_n(\sqrt{s} r) -
  \ds\frac{I_n(\sqrt{s})}{K_n(\sqrt{s})}K_n(\sqrt{s}r) \right]\cos n\theta, & 1<r\leq R;\\[10pt]
\ds\sum_{n=0}^{\infty} A_n\left[ \frac{I_n(\sqrt{s}R)}{K_n(\sqrt{s}R)} -
  \frac{I_n(\sqrt{s})}{K_n(\sqrt{s})} \right]K_n(\sqrt{s}r) \cos n\theta, & r>R,
\end{array}
\right.
\end{equation}
with constants $A_n$ determined from incorporation of the Dirac source to be
\begin{equation}
A_n = K_n(\sqrt{s}R) \left[  \int_{\theta=0}^{2\pi} \cos^2 (n\theta) \,d\theta \right]^{-1} = \left\{  \begin{array}{lc} \frac{1}{2\pi} K_0(\sqrt{s}R) & n = 0\\[5pt] \frac{1}{\pi} K_n(\sqrt{s}R) & n \geq 1 \end{array}
\right.
\end{equation}
Correspondingly, the flux over $r=1$ is given in series form by
\begin{equation}\label{eqn:fkuxLT}
 \hat{p}_r|_{r=1} = \frac{1}{2\pi} \frac{K_0(R\sqrt{s})}{K_{0}(\sqrt{s})} + \frac{1}{\pi}\sum_{n=1}^{\infty}\frac{K_n(R\sqrt{s})}{K_{n}(\sqrt{s})} \cos n\theta.
\end{equation}
To invert the Laplace transform of $\hat{p}_r|_{r=1}$, we must evaluate
the Bromwich integrals
\begin{align*}
\frac{1}{2\pi i} \int_{c- i \infty}^{c + i\infty}
  \frac{K_n(\sqrt{s}R)}{K_n(\sqrt{s})} e^{st} \, ds, \qquad n = 0,1,2,\ldots
\end{align*}
where $c$ is chosen to lie to the right of any poles of the integrand. Since the only singularity is a branch cut on the negative real axis, we deform the contour to a hairpin along the negative real axis and introduce the substitution $s = - w^2$. The integral becomes
\begin{align*}
\frac{1}{\pi i}\int_{0}^{\infty}\left[\frac{K_n(i\omega
R)}{K_n(i\omega)}-\frac{K_n(-i\omega R)}{K_n(-i\omega)}\right]\omega
e^{-\omega^2t}\,d\omega =
\frac{2}{\pi}\int_{0}^{\infty}\left[\frac{J_n(\omega)Y_n(\omega
R)-J_n(\omega
R)Y_n(\omega)}{Y_{n}^{2}(\omega)+J_{n}^{2}(\omega)}\right]\omega
e^{-\omega^2t}\,d\omega,
\end{align*}
where we have used $K_n(-iz) = \frac{\pi}{2} [-Y_n(z) + i J_n(z)]$. The expression for the flux $\mathcal{J}(t,\theta) = p_r|_{r=1}$ is now
\bsub
\begin{equation}
\mathcal{J}(t,\theta) = \frac{1}{2\pi} \chi_0(t)  + \frac{1}{\pi}\sum_{n=1}^{\infty}\chi_n(t) \cos n\theta, \qquad \theta \in (-\pi,\pi), \quad t>0,
\end{equation}
where the coefficients are
\begin{equation}
\chi_n(t) =
  \frac{2}{\pi}\int_{0}^{\infty}\left[\frac{J_n(\omega)Y_n(\omega
  \mu)-J_n(\omega
  \mu)Y_n(\omega)}{Y_{n}^{2}(\omega)+J_{n}^{2}(\omega)}\right]\omega
  e^{-\omega^2 t}\,d\omega.
\end{equation}
\esub
The total flux to the inner disk, and distribution of arrival times, is given by
\begin{equation}\label{eqn:timepdf}
\chi_0(t) = \int_0^{2\pi} \mathcal{J}(t,\theta) \, d\theta =
  \frac{2}{\pi} \int_{0}^{\infty}\left[\frac{J_0(\omega)Y_0(\omega
  R)-J_0(\omega
  R)Y_0(\omega)}{Y_{0}^{2}(\omega)+J_{0}^{2}(\omega)}\right]\omega
  e^{-\omega^2t}\,d\omega.
\end{equation}
Returning to dimensional time through \eqref{eqn:rescaled_vars}, we have that $S(t) = \frac{D}{a^2} \chi_0(\frac{D}{a^2} t)$. To determine $C(t)$ and $P(t) = 1-C(t)$, we note from $P'(t) = -S(t)$ that
\begin{equation}
P(t) = 1 - \int_0^{\frac{D}{a^2} t} \chi_0(\eta)\,d\eta, \qquad C(t) =
  \int_0^{\frac{D}{a^2} t} \chi_0(\eta)\,d\eta.
\end{equation}

For the dimensional arrival time $t = t^{\ast}$, the conditional distribution of arrival angles $\theta$ is then given by $\mathcal{J}(\frac{D}{a^2}t^{\ast},\theta)/\chi_0(\frac{D}{a^2}t^{\ast})$ with cumulative distribution
\begin{equation}
F(\theta;t_{\ast}) = \int_{-\pi}^{\theta}
  \frac{\mathcal{J}(\frac{D}{a^2}t_{\ast},\eta)}{\chi_0(\frac{D}{a^2}t_{\ast})}\,d\eta = \frac{\theta+\pi}{2\pi} + \frac{1}{\pi}\sum_{n=1}^{\infty} \frac{\chi_n(\frac{D}{a^2}t_{\ast})}{n\chi_0(\frac{D}{a^2}t_{\ast})} \sin n \theta, \qquad \theta\in(-\pi,\pi).
\end{equation}

\section*{Acknowledgments}
A.~E.~Lindsay was supported by NSF grant DMS-1815216. B.~Quaife was
supported by NSF grant DMS-2012560.

\bibliographystyle{siam}
%\bibliography{newbib}

\end{document}